\documentclass{amsart}
\numberwithin{equation}{section}
\usepackage[utf8]{inputenc}
\usepackage[T1]{fontenc} % this and the next package needed to use ogoneks
\usepackage{lmodern} % TODO: does this work with other fonts?
\usepackage{geometry}
\usepackage[usenames, dvipsnames]{xcolor}
\usepackage{graphicx}
\usepackage{subcaption}
\usepackage{mathtools}
\usepackage{amssymb}
\usepackage{amsthm}
\usepackage{xparse}
\usepackage{xspace}
\usepackage{mathrsfs}

\usepackage{setspace}
\usepackage[backref=page]{hyperref}
\usepackage[all]{hypcap}
\usepackage[capitalize, noabbrev]{cleveref}
\usepackage{booktabs}

\usepackage{adamsmacros}
\newcommand{\C}{\mathbb C} % complex numbers
\AtBeginDocument{ % some fonts redefine this (e.g. charter)

}
\newcommand{\N}{\mathbb N} % natural numbers
\newcommand{\R}{\mathbb R} % real numbers
\newcommand{\Z}{\mathbb Z} % integers
\newcommand{\cA}{\mathcal A} % Steenrod algebra, etc.

\newcommand{\GL}{\mathrm{GL}} % general linear
\newcommand{\SL}{\mathrm{SL}} % special linear
\AtBeginDocument{ % redefined by mathdesign
 \renewcommand{\O}{\mathrm O} % orthogonal
}
\newcommand{\SO}{\mathrm{SO}} % special orthogonal
\newcommand{\Spin}{\mathrm{Spin}} % spin
\newcommand{\Pin}{\mathrm{Pin}} % pin
\newcommand{\Spinc}{\relax\ifmmode{\Spin^c}\else Spin\textsuperscript{$c$}\xspace\fi}

\newcommand{\Pinc}{\relax\ifmmode{\Pin^c}\else Pin\textsuperscript{$c$}\xspace\fi}

\newcommand{\Pinp}{\relax\ifmmode{\Pin^+}\else Pin\textsuperscript{$+$}\xspace\fi}
\newcommand{\pinp}{pin\textsuperscript{$+$}\xspace}
\newcommand{\Pinm}{\relax\ifmmode{\Pin^-}\else Pin\textsuperscript{$-$}\xspace\fi}
\newcommand{\pinm}{pin\textsuperscript{$-$}\xspace}

\makeatletter
\usepackage[mathscr]{eucal}
\def\instring#1#2{TT\fi\begingroup
  \edef\x{\endgroup\noexpand\in@{#1}{#2}}\x\ifin@}
\def\isuppercase#1{%
  \instring{#1}{ABCDEFGHIJKLMNOPQRSTUVWXYZ}%
}%
\newcommand{\C@lIfUpper}[1]{
 \if\isuppercase{#1}\mathscr{#1}%
 \else #1%
 \fi
}
\newcommand{\cat}[1]{\mathit{\@tfor\next:=#1\do{\C@lIfUpper{\next}}}}
\makeatother
\newcommand{\newcat}[2]{\newcommand{#1}{\cat{#2}}}
\AtBeginDocument{
 \newcat{\fC}{C}
 \newcat{\fD}{D}
 \newcat{\Set}{Set} % sets
 \newcat{\Grp}{Grp} % groups
 \newcat{\Gpd}{Gpd} % groupoids
 \newcat{\Ab}{Ab} % abelian groups
 \newcat{\Ring}{Ring} % rings
 \newcat{\Mod}{Mod} % modules (over a ring)
 \newcat{\Alg}{Alg} % algebras (over a ring)
 \newcat{\Vect}{Vect} % vector spaces (over a field)
 \newcat{\sVect}{sVect} % super-vector-spaces
 \def\Top{\cat{Top}} % topological space (sometimes already defined, e.g. by kpfonts)
}
\newtheorem{thm}[equation]{Theorem}
\newtheorem{lem}[equation]{Lemma}
\newtheorem{cor}[equation]{Corollary}
\newtheorem{prop}[equation]{Proposition}

\theoremstyle{definition}
\newtheorem{exm}[equation]{Example}
\newtheorem{defn}[equation]{Definition}

\theoremstyle{remark}
\newtheorem{rem}[equation]{Remark}
\AddToHook{env/thm/begin}{\crefalias{equation}{thm}}
\AddToHook{env/lem/begin}{\crefalias{equation}{lem}}
\AddToHook{env/prop/begin}{\crefalias{equation}{prop}}
\AddToHook{env/cor/begin}{\crefalias{equation}{cor}}
\AddToHook{env/exm/begin}{\crefalias{equation}{exm}}
\AddToHook{env/rem/begin}{\crefalias{equation}{rem}}
\AddToHook{env/defn/begin}{\crefalias{equation}{defn}}

\crefname{thm}{Theorem}{Theorems}
\crefname{lem}{Lemma}{Lemmas}
\crefname{cor}{Corollary}{Corollaries}
\crefname{prop}{Proposition}{Propositions}
\crefname{ex}{Exercise}{Exercises}
\crefname{exmx}{Example}{Examples}
\crefname{defn}{Definition}{Definitions}
\crefname{claim}{Claim}{Claims}
\crefname{remx}{Remark}{Remarks}
\crefname{goal}{Goal}{Goals}
\crefname{ansatz}{Ansatz}{Ansatzes}
\crefname{notn}{Notation}{Notation}

\crefformat{section}{§#2#1#3}
\newcommand{\nocontentsline}[3]{}
\newcommand\stoptoc{%
   \let\origcontentsline\addcontentsline
   \let\addcontentsline\nocontentsline
}
\newcommand\resumetoc{%
   \let\addcontentsline\origcontentsline
}

\newcommand{\term}{\emph} % e.g. "The \term{trace} is defined to be..."

\newcommand{\RenewMathOperator}[2]{\renewcommand{#1}{\operatorname{#2}}}
\let\shortmapsto\mapsto
\renewcommand{\mapsto}{\mathchoice{\longmapsto}{\shortmapsto}{\shortmapsto}{\shortmapsto}}
\AtBeginDocument{

}
\DeclarePairedDelimiter\paren{(}{)}
\DeclarePairedDelimiter\ang{\langle}{\rangle}
\DeclarePairedDelimiter\abs{\lvert}{\rvert}

\DeclarePairedDelimiter\bkt{[}{]}
\DeclarePairedDelimiter\set{\{}{\}}

\makeatletter
\let\oldparen\paren
\def\paren{\@ifstar{\oldparen}{\oldparen*}}
\let\oldbkt\bkt
\def\bkt{\@ifstar{\oldbkt}{\oldbkt*}}
\makeatother
\newcommand{\newoperator}[1]{\expandafter\DeclareMathOperator\csname #1\endcsname{\operatorname{#1}}}
\newoperator{Aut}
\newoperator{coker}
\newoperator{End}
\newoperator{Ext}
\newoperator{Hom}
\AtBeginDocument{ % some packages redefine this
 \RenewMathOperator{\Im}{Im} % also complex analysis
}
\newoperator{Pic}
\newoperator{rank}
\newoperator{Sym}
\newoperator{Tor}
\newoperator{tr}
\newcommand{\bl}{\text{--}}

\newcommand{\Cl}{\mathit{C\ell}}

\newcommand{\KO}{\mathit{KO}}

\newcommand{\ko}{\mathit{ko}}

\newcommand{\Sq}{\mathrm{Sq}}

\newcommand{\NewThomSpectrum}[1]{\expandafter\newcommand\csname M#1\endcsname{M\mathrm{#1}}}
\newcommand{\NewMTSpectrum}[1]{\expandafter\newcommand\csname MT#1\endcsname{MT\mathrm{#1}}}
\newcommand{\BothThomSpectra}[1]{\NewThomSpectrum{#1}\NewMTSpectrum{#1}}
\BothThomSpectra{O}
\BothThomSpectra{SO}
\BothThomSpectra{Spin}
\BothThomSpectra{Pin}
\BothThomSpectra{U} % TODO: I haven't tested these
\AtBeginDocument{
 \RenewMathOperator{\Re}{Re}
}
\AtBeginDocument{
 \usepackage{microtype}
}
\usepackage{hyperref}
\hypersetup{
 colorlinks,
 linkcolor={red!50!black},
 citecolor={green!50!black},
 urlcolor={blue!80!black}
}

\usepackage{tikz-cd}

\newcommand{\NewCommenter}[2]{\expandafter\newcommand\csname #1\endcsname[1]{\textcolor{#2}{[{\bfseries #1:} ##1]}}}

\NewCommenter{Arun}{BrickRed}
\NewCommenter{CK}{OrangeRed}
\NewCommenter{Nat}{MidnightBlue}
\NewCommenter{Luuk}{Fuchsia}

\newcommand{\ftens}{\mathbin{\hat\times}}

\newcat{\sAlg}{sAlg}
\newcat{\sLine}{sLine} % added by CK

\definecolor{softred}{rgb}{0.92, 0, 0.17}

\newtheorem{notn}[equation]{Notation}

\DeclareDocumentCommand{\shortexact}{s O{} O{} mmmm}{
\IfBooleanTF{#1}{ % if star
\begin{tikzcd}[ampersand replacement=\&]
        {1} \& {#4} \& {#5} \& {#6} \& {1#7}
        \arrow[from=1-1, to=1-2]
        \arrow["#2", from=1-2, to=1-3]
        \arrow["#3", from=1-3, to=1-4]
        \arrow[from=1-4, to=1-5]
\end{tikzcd}
}{ % no star
\begin{tikzcd}[ampersand replacement=\&]
        {0} \& {#4} \& {#5} \& {#6} \& {0#7}
        \arrow[from=1-1, to=1-2]
        \arrow["#2", from=1-2, to=1-3]
        \arrow["#3", from=1-3, to=1-4]
        \arrow[from=1-4, to=1-5]
\end{tikzcd}
}}

\definecolor{ceruleanblue}{rgb}{0.16, 0.32, 0.75} % added by CK

\newcommand{\EMon}{\mathbb{M}}
\newcommand{\StMod}{{\cat{StMod}}}

\title{A tale of two periodicities}

\author{Arun Debray}
\address{Department of Mathematics, The University of Kentucky, 719 Patterson Office Tower,
Lexington, \indent KY 40506, USA}
\email{\href{mailto:a.debray@uky.edu}{a.debray@uky.edu}}
\urladdr{\href{https://adebray.github.io/}{https://adebray.github.io/}}

\author{Cameron Krulewski}
\address{Perimeter Institute for Theoretical Physics,
31 Caroline St North,
Waterloo, ON,
N2L 2Y5, Canada}
\email{\href{ckrulewski@perimeterinstitute.ca}{ckrulewski@perimeterinstitute.ca}}
\urladdr{\href{https://cakrulewski.github.io/}{https://cakrulewski.github.io/}}
\author{Natalia Pacheco-Tallaj}
\address{Department of Mathematics,
University of California, Berkeley,
970 Evans Hall,
Berkeley, CA 94720, \indent USA}
\email{\href{mailto:nmpt@berkeley.edu}{nmpt@berkeley.edu}}
\urladdr{\href{https://natipt.github.io/}{https://natipt.github.io/}}

\subjclass[2020]{\href{https://mathscinet.ams.org/msc/msc2020.html?t=19L41&btn=Current}{19L41}, \href{https://mathscinet.ams.org/msc/msc2020.html?t=55S10}{55S10}}
\date{\today}

\begin{document}
\maketitle

\begin{abstract}
We provide a common explanation for two $4$-periodic phenomena arising in the study of connective real $K$-theory: first, the family of twists of $\ko$ over $B\Z/2$, and second, the torsion in the Picard group of $Q_1$-local $\cA(1)$-modules. These phenomena have been previously studied, as they lead to $4$-periodic patterns in Smith homomorphisms, resp.\ in the $\ko$-homology of $(B\Z/2)^n$. In this paper, we exhibit both of these $4$-periodicities as consequences of an isomorphism of two extraspecial groups of order $32$.
\end{abstract}

\tableofcontents

%%%%%%%%%%%%%%%%%%%%%
%%%%%%%%%%%%%%%%%%%%%

\setcounter{section}{-1}
\section{Introduction}

In topology, many occurrences of the number $8$ have the same origin: Bott periodicity in $\KO$-theory. Quite different-looking $8$s, to do with quadratic forms, bordism groups, or Clifford algebras, can be identified by theorems in algebra, geometry, or topology~\cite{Ranicki_8}. %\Arun{Transition between $8$ and $4$ kind of abrupt, maybe we can say something re: the connection between all these $8$s inspired us to find a connection between these $4$s?}

In prior work joint with Devalapurkar, Liu, Stehouwer, and Thorngren~\cite{debray_long_2024, debray_smith_2024, BottSpiral}, we learned of two different \emph{4}-periodic phenomena related to connective $\ko$-theory. As a parallel to the connections between all these $8$s, we wondered whether these two $4$-periodicities can be identified.

These two periodicities are:
\begin{itemize}
    \item the $4$-periodic family of twists of spin bordism and $\ko$-theory over $B\Z/2$, as discussed in~\cite{kapustin_fermionic_2015,TY19,HKT19,COSY,WWZ19,debray_long_2024,debray_smith_2024}, and
    \item the $4$-periodic pattern in the $\ko$-homology of $B(\Z/2)^n$ arising from the Picard group of $Q_1$-local stably invertible $\cA(1)$-modules, as studied by~\cite{Yu95, Bru14}.
\end{itemize}
The former periodicity comes from the Thom isomorphism: both $\ko$-theory and spin bordism have Thom isomorphisms for \textit{spin} vector bundles, but either theory evaluated on the Thom space of a non-spin vector bundle can be identified with a \emph{twisted} bordism or $\ko$-group over the base space.
Over $B\Z/2$, there are four inequivalent twists, corresponding to the trivial rank-zero bundle, $\sigma-1$, $2\sigma-2$, and $3\sigma-3$, where $\sigma\to B\Z/2$ is the universal bundle. Once we take four copies of $\sigma$, we again have an untwisted Thom isomorphism: $4\sigma$ admits a spin structure.
In quantum physics, this collection of four twisted structures corresponds to a set of four types of fermionic systems: those with a certain symmetry operator that is either unitary or antiunitary and which squares either to 1 or to $-1$; see e.g.\ the references above for details.

Meanwhile, one can discover the latter periodicity when computing $\ko$-homology via the Adams spectral sequence. The input data to this computation is a module over the subalgebra $\cA(1) = \ang{\Sq^1, \Sq^2}$ of the mod $2$ Steenrod algebra. Yu~\cite{Yu95} showed that, after localizing at the collection of $Q_1$-Margolis homology isomorphisms (see \cref{marghom}), the category of $\cA(1)$-modules has exactly four isomorphism classes of stably invertible objects up to suspension.
%after
%$Q_1$-localizing, there
%are exactly four stably invertible, $\cA(1)$-modules in the Picard group 
%\CK{hopefully said that correctly} \Arun{I tried to rephrase it, what do you think?}. \CK{need to recall that $Q_1$-localizing is a common technique?} \Arun{I think it's obscure enough that we should either say something or (more likely) forward-ref the defn}

A priori, these two $4$-periodic phenomena do not have much in common. We began to wonder whether they were related in prior work of two of us and Stehouwer in~\cite{BottSpiral}, where both phenomena played a role in computing twisted $\ko$-theory of $(B\Z/2)^n$ and a corresponding summand of twisted spin bordism. 

%In particular, we sought a common explanation for these two $4$-periodic phenomena.
%This latter periodicity clearly informs the former: in \cite{BottSpiral}, the first two authors used the periodicity of these $Q_1$-local $\cA(1)$-modules as an input to an Adams spectral sequence computation of the .
%\CK{say something about the $\times 16$ here?} \Arun{what sort of thing did you have in mind?}
%However, we suspected that the two periodicities shared a deeper origin: one relating to a certain class of extraspecial 2-groups.

In this paper, we identify the common origin of both periodicities.
Specifically, we 
define a commutative monoid $\EMon$ whose elements are certain extraspecial $2$-groups and whose monoid operation, the \term{fermionic tensor product}, is a generalization of the central product (see \cref{ftens_char}). $\EMon$ maps to our two periodicities in the following ways.
\begin{itemize}
    \item Taking twisted cohomology defines a homomorphism from $\EMon$ to the Picard group of $Q_1$-local $\cA(1)$-modules.
    %A standard construction to extract an $\cA(1)$-module from TODO is group completion.
    \item Passing to the \term{diagonal subgroup} (see \cref{elk_exm}) defines a twist of $\ko$ over $B\Z/2$.
\end{itemize}
%We prove in \cref{period_compare} that 
% comprised of bigraded extraspecial super 2-groups (which we denote $E_{\ell,k}$ for $(\ell,k)\in \Z^2$) and 
We prove in \cref{period_compare} that the torsion $\Z/4$ subgroup of the group completion of $\EMon$ admits surjective maps to the corresponding $\Z/4$ subgroups in the $\Z/2$-twists of spin bordism and the stable $Q_1$-local Picard group, respectively.

Now we go into more detail about our main results.
In \cref{elk_exm}, we define extraspecial $2$-groups $E_{\ell,k}$ of order $2^{\ell+k+1}$ for all $\ell,k\ge 0$. Heuristically, $E_{\ell,k}$ is defined by taking $\ell$ copies of $\Pin_1^+$ and $k$ copies of $\Pin_1^-$, then asking them to commute in the graded sense. Then (\cref{EMon_defn}) we define $\EMon\coloneqq\{E_{\ell,k}\mid \ell, k\ge 0\}$. The fermionic groups $E_{\ell,k}$ within $\EMon$ correspond, after taking a twisted group algebra, to the Clifford algebras $\Cl_{\ell,k}$. There is an isomorphism $E_{4,0}\cong E_{0,4}$ (see \cref{seed}), corresponding to the isomorphism  $\Cl_{0,4}\cong \Cl_{4,0}$. This yields an important relation in $\EMon$, which leads to our first main result.
%Define $\EMon$ here
%We obtain $\EMon$ by taking the set of $E_{\ell,k}$ for all $(\ell,k)$ and then quotienting by the relation corresponding to the isomorphism
%
%\begin{itemize}
%    \item define $\EMon$?
%    \item define $E_{\ell,k}$?
%\end{itemize}%
%
%\CK{and then on to the summary? save symmetry group interpretation for section 1? save mention of Bott spiral for after the summary or a remark later in the paper?}
%
%\noindent\hrulefill

\begin{thm}
\label{Q1mainthm}
The monoid homomorphism
\begin{equation}
    \EMon \xrightarrow[\eqref{our_phi_defn}]{\Phi} \Pic(\StMod_{\cA(1)}[Q_1^{-1}])\xrightarrow[\eqref{yu_thm}]{y} \Z\oplus\Z/4
\end{equation}
defined in \cref{our_phi_defn} is the group completion.
\end{thm}
The map $\Phi$ is the twisted cohomology map mentioned above.
\Cref{yu_thm} includes that $y$ is an isomorphism, so $\Phi$ on its own is also a group completion; we compose with $y$ to make the $4$-periodicity easier to see.

The $\Z/4$ torsion summand of $\Pic(\StMod_{\cA(1)}[Q_1^{-1}])$ constructed in Yu's theorem~\cite{Yu95} (\cref{yu_thm}) is somewhat mysterious, so it is nice that
\cref{Q1mainthm} gives a group-theoretical explanation for it.

%\TODO: words here re: hey it's cool we found an origin for the $\Z/4$ in the $Q_1$-local Picard group

In \cref{monoid_hom_cor}, we show that the map sending $E_{\ell,k}$ to the twisting data of the diagonal subgroup $\Delta_{\ell,k}\subset E_{\ell,k}$ defines a monoid homomorphism
\begin{equation}
    \delta\colon\EMon\longrightarrow [B\Z/2, B\O/B\Spin].
\end{equation}
Since $[B\Z/2, B\O/B\Spin]$ is an abelian group, by \cref{Q1mainthm} there is a unique group homomorphism 
\begin{equation}
\label{tilde_delta}
\widetilde\delta\colon \Z\oplus \Z/4\longrightarrow [B\Z/2, B\O/B\Spin]
\end{equation}
satisfying $\delta = \widetilde\delta\circ y \circ \Phi$. This suggests an obvious way to compare these two periodicities: if $t\colon \Z/4\hookrightarrow\Z\oplus\Z/4$ is the inclusion of the torsion subgroup, we could take the comparison map
\begin{equation}
\label{naive_comparison}
    \Psi\colon \Z/4\overset t\hookrightarrow \Z\oplus\Z/4 \xrightarrow{~\widetilde\delta~} [B\Z/2, B\O/B\Spin] \xrightarrow[\eqref{second_Z4}]{\cong} \Z/4.
\end{equation}
However, we show in \cref{calc_psi} that this map is \emph{not} an isomorphism! Instead, it sends $1\mapsto 2$.

Nonetheless, these two $4$-periodicities are still related, as the fermionic group $E_{1,0}$ maps to generators of both $\Z/4$s. In \S\ref{s:diag} we comment on why the comparison map $\Psi$ does not work: there is an ``off-by-one error'' in that $\Phi(E_{1,0})$ is a suspension of the generator of the torsion subgroup of $\Pic(\StMod_{\cA(1)}[Q_1^{-1}])$, and this is the obstruction to $\Psi$ being an isomorphism.

Once we correct for this shift, the suspense is over, and we obtain a correspondence between these two periodicities.
% showing they are related. 
Specifically, in \cref{desuspend}, we construct a projection homomorphism $\varsigma$ from $\Pic(\StMod_{\cA(1)}[Q_1^{-1}])$ onto its torsion subgroup by suspending or desuspending a $Q_1$-locally invertible $\cA(1)$-module until its image in the Picard group is torsion. This suffices to fix the ``off-by-one'' issue noted above.

%\TODO: below, the domain of $\widetilde\delta$ is $\Z\oplus\Z/4$, not the Picard group. So maybe precompose with $y^{-1}$ or something
\begin{thm}
\label{period_compare}
%Let $t\colon\Z/4\hookrightarrow \Z\oplus\Z/4$ be the inclusion of the torsion subgroup. Then the following composition is an isomorphism:
%
%
Both maps $\widetilde\delta\circ y$ and $\varsigma$ in the diagram
\begin{equation}
\label{weird_span}
\begin{tikzcd}[column sep=-3em]
% https://q.uiver.app/#q=WzAsNixbMSwwLCJcXG1hdGhiYiBMIl0sWzAsMSwiW0JcXFovMiwgQlxcTy9CXFxTcGluXSJdLFswLDIsIlxcWi80Il0sWzIsMSwiXFxQaWMoXFxjYXR7U3RNb2R9X3tcXGNBKDEpfVtRXzFeey0xfV0pIl0sWzIsMywiXFxaLzQiXSxbMiwyLCJcXFBpYyhcXGNhdHtTdE1vZH1fe1xcY0EoMSl9W1FfMV57LTF9XSlfe1xcbWF0aHJte3RvcnN9fSJdLFswLDEsIlxcZGVsdGF8X3tcXG1hdGhiYiBMfSIsMl0sWzEsMiwiXFxjb25nIiwyXSxbMCwzLCJcXFBoaXxfe1xcbWF0aGJiIEx9Il0sWzMsNSwiXFx2YXJzaWdtYSJdLFs1LDQsIlxcY29uZyJdXQ==
	& {\Pic(\StMod_{\cA(1)}[Q_1^{-1}])} & \\
	{[B\Z/2, B\O/B\Spin]} && {\Pic(\StMod_{\cA(1)}[Q_1^{-1}])_{\mathrm{tors}}} \\
	{\Z/4} && {\Z/4}
	\arrow["{\widetilde\delta\circ y}"', "{(\ref{tilde_delta},\,\ref{yu_thm})}"{pos=0.45}, from=1-2, to=2-1]
	\arrow["\varsigma", "\eqref{desuspend}"', from=1-2, to=2-3]
	\arrow["\cong"', "\eqref{second_Z4}", from=2-1, to=3-1]
	\arrow["\cong", "\eqref{yu_thm}"', from=2-3, to=3-3]
\end{tikzcd}\end{equation}
send $\Phi(E_{1,0})$ to a generator, hence are surjective.
\end{thm}
Thus, after applying $\varsigma$ to desuspend $\Phi(E_{1,0})$, restricting from $E_{\ell,k}$ to its diagonal subgroup $\Delta_{\ell,k}$---i.e.\ applying the map $\widetilde\delta\circ y$ in~\eqref{weird_span}---identifies the two $\Z/4$s we study in this note, hence also the two $4$-periodicities they define.
%In other words, restricting from $E_{\ell,k}$ to its diagonal subgroup $\Delta_{\ell,k}$ defines an isomorphism between the two $\Z/4$s we study in this note, hence between the two $4$-periodicities that they define. 
As briefly discussed in~\cite[Remark 5.46]{BottSpiral}, our primary goal in this note is to understand whether these two four-periodic phenomena were related, and \cref{period_compare} offers a strong answer: in this sense, they're the \emph{same} $4$-periodicity.  

It would be interesting to upgrade the correspondence in \cref{period_compare} into a canonically defined group isomorphism $\Z/4\xrightarrow\cong\Z/4$---one can use \cref{period_compare} to produce an isomorphism like this, but only after making choices.

% \CK{maybe sneak in something \url{https://poets.org/teach-poem/excerpt-tale-two-cities} ``It was the best of times, it was the worst of times, it was the age of wisdom, it was the age of foolishness, it was the epoch of belief, it was the epoch of incredulity, it was the season of Light, it was the season of Darkness, it was the spring of hope, it was the winter of despair, we had everything before us, we had nothing before us, we were all going direct to Heaven, we were all going direct the other way---in short,the period was so far like the present period, that some of its noisiest authorities insisted on its being received, for good or for evil, in the superlative degree of comparison only''...}

\begin{rem}
The four periodicities in this paper yield physical interpretations.
Just as Bott periodicity informs the classification of topological insulators and superconductors~\cite{kitaev_periodic_2009}, the four periodicities in this paper give rise to a periodic structure in the classification of more general, interacting versions of these fermionic quantum systems. This structure was originally observed by Queiroz--Khalaf--Stern~\cite{queiroz_dimensional_2016} and two of us, joint with Stehouwer, studied it mathematically in~\cite{BottSpiral}.
%In particular,
%    The groups $E_{\ell,k}$ came to our attention in that work: they are the symmetry groups of a collection of physical systems, collectively called the \term{Bott spiral}, studied by ~\cite{queiroz_dimensional_2016}. See~\cite{BottSpiral} for more information.
\end{rem}

\subsection*{Outline}
We begin in \S\ref{s:fermionic} with a brief introduction to fermionic groups. Though there are already good introductions to fermionic groups in~\cite{Ben88, Sto98, stehouwer_interacting_2022, BottSpiral}, in this paper we need a symmetric monoidal $\infty$-category of discrete fermionic groups, so take a more homotopical viewpoint than those references. We then introduce our fundamental examples $E_{\ell,k}$ (\cref{elk_exm}), use them to build the commutative monoid $\EMon$ (\cref{EMon_defn}) and the diagonal subgroups $\Delta_{\ell,k}$~\eqref{defndiag}, and show that taking the twisting data of the diagonal defines a monoid homomorphism $\EMon\to[B\Z/2, B\O/B\Spin]$ (\cref{monoid_hom_cor}).

In \S\ref{ss:first_Z4}, we introduce $Q_k$-local stable module categories of $\cA(n)$-modules, then give the first of our two periodicities, which arises as the $\Z/4$ subgroup of the $Q_1$-local stable Picard group of $\cA(1)$-modules. The main result in this subsection is Yu's calculation~\cite[Theorem 2.1]{Yu95} (here \cref{yu_thm}) that this Picard group is isomorphic to $\Z\oplus\Z/4$.

In \S\ref{ss:second_Z4}, we compute $[B\Z/2, B\O/B\Spin]\cong\Z/4$ (\cref{second_Z4}) and interpret it as giving a $4$-periodic sequence of twists of real $K$-theory, hence a $\Z/8\times\Z/4$-grading on the twisted $\KO$-cohomology ring of $B\Z/2$.

In \S\ref{s:assemble}, we construct the monoid homomorphism $\Phi\colon\EMon\to\Pic(\StMod_{\cA(1)}[Q_1^{-1}])$ that links these two periodicities together.

Finally, in \S\ref{s:proofs}, we prove our main theorems relating these two $4$-periodicities: \cref{Q1mainthm,period_compare,calc_psi}.

\subsection*{Acknowledgments}
This paper was born as a spinoff of~\cite{BottSpiral}, joint between two of us and Luuk Stehouwer. We would especially like to thank Luuk Stehouwer for many helpful discussions regarding the groups $E_{\ell,k}$ and their role in these two periodicities.
%for many helpful discussions on the content of this paper w
%\TODO: Luuk for previous collaboration, funding, who else? 
We also want to thank Thorger Geiß, Theo Johnson-Freyd, and Noam Nissan for helpful conversations on other aspects of this paper.

CK was supported by the National Science Foundation under Award No. 2602343, the Simons Foundation, and the Killam Trusts during the writing of this paper. CK is based at the Perimeter Institute for Theoretical Physics;
research at Perimeter is supported by the Government of Canada through Industry Canada and by the Province of Ontario through the Ministry of Research \& Innovation.
Any opinions, findings, and conclusions or recommendations expressed in the materials in this work are those of the authors and do not necessarily reflect the views of the National Science Foundation.
NPT was supported by the UC President's Postdoctoral Fellowship.

\subsubsection*{Usage of AI}
One of us consulted Claude during literature review, but we did not use AI to write or prove the results in this paper.
%
%\TODO
%\CK{I hope we can say something like %``''} \Arun{I did not use it at all}

%%%%%%%%%%%%%%%%%%%%%
%%%%%%%%%%%%%%%%%%%%%

\section{Fermionic groups}
\label{s:fermionic}
%[\TODO: flavor text here! probably say ``fermionic'' comes from a physics application; see our other paper(s). but it's also useful for precisely tracking products of twists of $K$-theory and spin bordism, which is useful in index theory, etc]
%\CK{I put in \cref{ferm_grp_name_why}. maybe that is sufficient?}
%
%
Let $A$ be an $E_\infty$-algebra and let $\Top_{A}$ denote the slice $\infty$-category of spaces $X$ with a map $f_X\colon X\to A$. Maps in this $\infty$-category are continuous maps of spaces that commute with the maps $f_X$ to $A$.

\begin{prop}[{Ando--Blumberg--Gepner~\cite[Proposition 6.12]{ABG18}}]
\label{slice_tensor}
For $A$ an $E_\infty$-algebra in spaces, there is a natural symmetric monoidal $\infty$-category structure on $\Top_{/A}$ whose value on two objects $(X, f_X\colon X\to A)$ and $(Y, f_Y\colon Y\to A)$ is the Cartesian product and external sum
\begin{equation}
    (X\times Y, \mu\circ (f_X\times f_Y)\colon X\times Y\to A),
\end{equation}
where $\mu\colon A\times A\to A$ is the product on $A$.
\end{prop}

See
Beardsley--Luecke--Morava~\cite[Remark 5.10]{beardsley2023brauer} for a little more detail on the relationship between Ando--Blumberg--Gepner's construction and $\Top_{/A}$. Naturality in \cref{slice_tensor} is with respect to composition with maps of $E_\infty$-spaces.

\begin{defn}
\label{bo_bspin}
Give $B\O$ and $B\Spin$ the $E_\infty$-structures arising from the direct sum of vector bundles, resp.\ of spin vector bundles. Then, let $B\O/B\Spin$ be the cofiber, in $E_\infty$-spaces, of the forgetful map $B\Spin\to B\O$.
\end{defn}

\begin{lem}[{\cite[Proposition 1.37 and Lemma 1.41]{DY23}}]
\label{bobsp_split}
As spaces, there is a preferred homotopy equivalence $\phi\colon B\O/B\Spin\xrightarrow{\simeq} K(\Z/2, 1)\times K(\Z/2, 2)$, but this does not extend to an equivalence of $E_1$-spaces.
\end{lem}

\begin{defn}
\label{FGdefn}
The category $\cat{FermGrp}$ of \term{(discrete)\footnote{We emphasize ``discrete'' because, in many cases in prior literature, one works with a category of fermionic topological groups.} fermionic groups}
is the full sub-$\infty$-category of $\Top_{/(B\O/B\Spin)}$ on the objects $(X, f_X\colon X\to B\O/B\Spin)$ such that $X$ is connected and $1$-truncated.
\end{defn}
%[\TODO: footnote on the word ``discrete''?]

The property of being connected and $1$-truncated is preserved by Cartesian products, so the symmetric monoidal structure on $\Top_{/(B\O/B\Spin)}$ restricts to a symmetric monoidal $\infty$-category structure on $\cat{FermGrp}$.

\begin{rem}\label{ferm_grp_name_why}
    The name \textit{fermionic group} comes from physics, where these objects may be used to describe the symmetries of quantum systems of fermions; e.g. electrons.
\end{rem}

Now that we have a well-behaved category of fermionic groups, let's see how to actually work with them.
\begin{defn}
\label{fermgp}
Let $X\to B\O/B\Spin$ be a fermionic group. We define the following invariants associated to $X$:
\begin{enumerate}
    \item The \term{bosonic group} $G_b\coloneqq \pi_1(X)$. Since $X$ is connected and $1$-truncated, $X\simeq BG_b$.
    \item\label{thetapart} The composition
    \begin{equation}
        \theta\colon BG_b \simeq X\to B\O/B\Spin \to B\O/B\SO\simeq B\O_1,
    \end{equation}
    induced by the forgetful map $B\Spin\to B\SO$, can be reinterpreted as a class in $H^1(BG_b;\Z/2)\cong \Hom(G_b, \Z/2)$, which we will also call $\theta$.
    \item The map $B\O_1\to B\O$ classifying the universal bundle induces a map $B\O/B\Spin\to B\SO/B\Spin\simeq K(\Z/2, 2)$. Proceeding as in the previous bullet point, we obtain a class $\omega\in H^2(BG_b;\Z/2)$. This is equivalent data to an isomorphism class of extension of $G_b$ by $\Z/2$, which we call $G_f$.
\end{enumerate}
The pair $(\theta, \omega)$ is called the \term{twisting data} of $X$.
\end{defn}

\begin{lem}\label{lem:Gf_vs_twisting_data}
The data $(G_b, \theta, \omega)$ determine a fermionic group up to isomorphism in the homotopy category of $\cat{FermGrp}$.
\end{lem}
\begin{proof}
As discussed above, $G_b$ determines $X$ up to homotopy, since $X\simeq BG_b$. The homotopy equivalence $B\O/B\Spin\simeq K(\Z/2, 1)\times K(\Z/2, 2)$ mentioned above sends the map $f_X\colon X\to B\O/B\Spin$ to the twisting data $(\theta, \omega)\colon BG_b\to K(\Z/2, 1)\times K(\Z/2, 2)$, hence determines $f_X$ up to homotopy equivalence. To finish, observe that in the homotopy category of $\Top_{/Z}$ for a space $Z$, an isomorphism is a homotopy equivalence of spaces that commutes, up to homotopy, with the maps down to $Z$.
%homotopy-equivalent spaces with homotopic maps to $Z$ are isomorphic.
\end{proof}

We will generally refer to a fermionic group either as the data $(G_b, \theta, \omega)$, or just as $G_f$ when the rest of the data is clear from context. The monoidal product of $G_f$ and $H_f$ in $\cat{FermGrp}$, using the symmetric monoidal structure described above, will be denoted $G_f\ftens H_f$.

\begin{prop}
\label{ftens_char}
Let $G_f$ and $H_f$ be fermionic groups with twisting data $(G_b, \theta_G, \omega_G)$ and $(H_b, \theta_H, \omega_H)$, respectively. Then $G_f\ftens H_f$ may be characterized, up to isomorphism in $\mathrm{Ho}(\cat{FermGrp})$, in either of the following ways.
\begin{enumerate}
    \item\label{whitney_tens} The fermionic group $K_f$ whose bosonic group is $K_b\coloneqq G_b\times H_b$ and whose twisting data is
    \begin{equation}
    \label{whit_twist}
        (\theta_K, \omega_K) = (\theta_G + \theta_H, \omega_G + \theta_G\theta_H + \omega_H).
    \end{equation}
    \item\label{other_tens} The \term{fermionic tensor product} of $G_f$ and $H_f$ as described by Benson~\cite[\S 7]{Ben88}: as a set, the direct product $G_f\times H_f$, with multiplication
    \begin{equation}
(g_1 \mathbin{\hat{\otimes}} h_1) (g_2 \mathbin{\hat{\otimes}} h_2) = 
\begin{cases}
(-1)^F g_1 g_2 \mathbin{\hat{\otimes}} h_1 h_2, & \text{ if } \theta(h_1) = 1 \text{ and } \theta(g_2) = 1
\\
    g_1 g_2 \mathbin{\hat{\otimes}} h_1 h_2,  & \text{ if } \theta(h_1) = 0 \text{ or } \theta(g_2) = 0.
\end{cases}
\end{equation}
\end{enumerate}
\end{prop}
\begin{proof}
Part~\eqref{other_tens} follows from~\eqref{whitney_tens} and~\cite[Lemma 1.12]{BottSpiral}, so we focus on part~\eqref{whitney_tens}. The monoidal product on $\cat{FermGrp}$ is the Cartesian product $BG_b\times BH_b$, which is homotopy equivalent to $B(G_b\times H_b)$, so we may take $K_b\coloneqq G_b\times H_b$. To compute the twisting data of $G_f\ftens H_f$, it suffices to understand how the $E_\infty$-structure on $B\O/B\Spin$ behaves under the (non-$E_1$) homotopy equivalence $B\O/B\Spin\simeq K(\Z/2, 1)\times K(\Z/2, 2)$, which sends a map $B\O/B\Spin$ to a pair of degree-$1$ and degree-$2$ $\Z/2$-cohomology classes. In~\cite[(1.42)]{DY23} (see also~\cite[(10.9)--(10.11)]{freed_reflection_2021}), this identification is shown to produce the formula~\eqref{whit_twist}.
%The formula for the twisting data then follows from the identification of the $E_\infty$-structure on $B\O/B\Spin$ under the (non-$E_1$) homotopy equivalence
%Using the homotopy equivalence $B\O/B\Spin\simeq K(\Z/2, 1)\times K(\Z/2, 2)$
%As described in, the identifications of $B\O/B\Spin\to K(\Z/2, 1)\times K(\Z/2, 2)$ and maps to $K(\Z/2, n)$ with cohomology classes
%carry the external sum of maps $(\theta_G, \omega_G)$ and $(\theta_H, \omega_H($ to $B\O/B\Spin$ to the 
\end{proof}

\begin{rem}
\label{H1H2_twist}
For a fixed finite group $G_b$, one can consider the set of isomorphism classes of fermionic groups $G_f$ with bosonic group isomorphic to $G_b$. By \cref{lem:Gf_vs_twisting_data}, we may equivalently consider the collection of different twisting data over $BG_b$. As a set, this is $[BG_b, B\O/B\Spin]\cong H^1(BG_b;\Z/2)\times H^2(BG_b;\Z/2)$; using the grouplike $E_\infty$-structure on $B\O/B\Spin$, this becomes an abelian group, a direct summand of the \term{Brauer group} of $BG_b$~\cite{donovan_graded_1970, beardsley2023brauer}. Explicitly, the product in this group structure on the set $H^1\times H^2$ has the formula~\eqref{whit_twist}.

\textbf{Warning:} %\CK{put this in black?}
we will write elements of $[X, B\O/B\Spin]$ as elements of $H^1(X;\Z/2)\times H^2(X;\Z/2)$. However, the group operation~\eqref{whit_twist} is \emph{not} the product of the operations on $H^1$ and $H^2$.
\end{rem}

\begin{exm}[{$\Pin_1^\pm$}]
\label{pin_exms}
Recall that $\Pin_1^\pm$ is the group of units in the Clifford algebra on one generator $e$, with $e^2 = \pm 1$. Thus $\Pin^+_1\cong \Z/2\times\Z/2$ and $\Pin^-_1 \cong \Z/4$. In both cases, quotienting by $\{\pm 1\}$ gives $\O_1\cong \Z/2$, with $e$ mapping to the generator. We thus obtain two fermionic groups $\Pin^\pm_1$ with bosonic part $(\Pin^\pm_1)_b \cong \Z/2$.

Write $H^*(B\Z/2;\Z/2) = \Z/2[x]$ with $\abs x =1$. 
\begin{itemize}
    \item We choose $\theta = x$ for both $\Pin_1^+$ and $\Pin_1^-$. This represents the fact that the quotient maps $\Pin_1^\pm\to\O_1$ are both surjective.
    \item The class $\omega$ is the extension class of $\Z/2\to \Pin^\pm_1\to \Z/2$. The extension is split for $\Pin^+_1$ and nonsplit for $\Pin^-_1$, so $\omega=0$ for $\Pin^+_1$ and $\omega = x^2$ for $\Pin^-_1$.
\end{itemize}
The twisting data for $\Pin^+_1$ is thus $(x,0)$ and for $\Pin^-_1$ is $(x,x^2)$. 
\end{exm}

\begin{exm}
The \term{Conway fermionic group} has $G_f = \mathrm{Co}_0$ and $G_b = \mathrm{Co}_1$, with $\theta$ trivial. Here $\mathrm{Co}_1$ is the sporadic simple group discovered by Conway in~\cite{conway_perfect_1968,conway_group_1969}, where it is called $\cdot 1$, and $\mathrm{Co}_0$ is its double cover, called $\cdot 0$ by Conway (\textit{ibid.}), which is the automorphism group of the Leech lattice. See Johnson-Freyd--Treumann~\cite[\S 6]{JFT20} and Albert--Okada--Tachikawa~\cite{albert_art_2026} for applications of this fermionic group to questions in physics.
\end{exm}

\begin{exm}[{\cite[Definition 1.36]{BottSpiral}}]
\label{elk_exm}
Define $E_{\ell,k}$ to be the fermionic tensor product
\begin{equation}
     E_{\ell,k}\coloneqq \underbracket{\Pin_1^+\ftens  \dotsb \ftens  \Pin_1^+}_{\text{$\ell$ copies}} \ftens 
     \underbracket{\Pin_1^-\ftens  \dotsb \ftens  \Pin_1^-}_{\text{$k$ copies}}.
\end{equation}
Thus $(E_{\ell,k})_b\cong (\Z/2)^{\ell+k}$, and
\begin{equation}
\label{add_mon}
E_{\ell,k}\ftens E_{\ell',k'}\cong E_{\ell+\ell', k+k'}.
\end{equation}
See~\cite[(1.40)]{BottSpiral} for the twisting data of $E_{\ell,k}$. Define the \term{diagonal subgroup} $\Delta_{\ell,k}\subset E_{\ell,k}$ to be the pullback
% https://q.uiver.app/#q=WzAsNCxbMCwwLCJcXERlbHRhX3tcXGVsbCxrfSJdLFsxLDAsIkVfe1xcZWxsLGt9Il0sWzEsMSwiKFxcWi8yKV57XFxlbGwra30iXSxbMCwxLCJcXFovMiJdLFszLDIsIjFcXG1hcHN0byAoMSxcXGRvdHNjLDEpIl0sWzAsM10sWzEsMl0sWzAsMV0sWzAsMiwiIiwxLHsic3R5bGUiOnsibmFtZSI6ImNvcm5lciJ9fV1d
\begin{equation}\label{defndiag}\begin{tikzcd}[column sep=4.5em]
	{\Delta_{\ell,k}} & {E_{\ell,k}} \\
	{\Z/2} & {(\Z/2)^{\ell+k}.}
	\arrow[from=1-1, to=1-2]
	\arrow[from=1-1, to=2-1]
	\arrow["\lrcorner"{anchor=center, pos=0.125}, draw=none, from=1-1, to=2-2]
	\arrow[from=1-2, to=2-2]
	\arrow["{1\mapsto (1,\dotsc,1)}", from=2-1, to=2-2]
\end{tikzcd}\end{equation}
$\Delta_{\ell,k}$ is a fermionic group with bosonic group $\Z/2$ and twisting data pulled back from that of $E_{\ell,k}\to (\Z/2)^{\ell+k}$.

The first examples are trivial: the inclusions $\Delta_{1,0}\to E_{1,0}$ and $\Delta_{0,1}\to E_{0,1}$ are isomorphisms. The inclusion $\Delta_{1,1}\to E_{1,1}$ is more interesting:
one can show that there is an isomorphism $g\colon D_8\to E_{1,1} = \Pin_1^+\ftens \Pin_1^-$ sending a $\pi/2$ rotation $r$ to $(1, e_1)$ and a reflection $s$ to $(e_1, 1)$. (Here we are using the description of the fermionic tensor product from part~\eqref{other_tens} of \cref{ftens_char}.)  $\Delta_{1,1}$ is the preimages of $(0,0)$ and $(1,1)$ under $E_{1,1}\to (\Z/2)^2$, so $\Delta_{1,1} = \{\pm (1,1), \pm (e_1, e_1)\}$.  Therefore $g^{-1}(\Delta_{1,1})$ is the order-$4$ group $\{1, r^2, s, sr^2\}$. In other words, $\Delta_{1,1}$ is ``trivial'': its bosonic group is $\Z/2$ and its twisting data is $(0, 0)$.
%[\TODO: maybe wordsmith a little]
%\CK{actually I think it's fine?}
%(Here, following Benson~\cite[\S 7]{Ben88}, we write the fermionic tensor product as the direct product of sets, equipped with the twisted multiplication
%\TODO: Introduce $E_{\ell,k}$, and call $\Delta_{\ell,k}$ the diagonal subgroup (whose nontrivial element is like $(1,1,...,1)\in (\Z/2)^{\ell+k}$). Explicitly work this out for $\Delta_{1,1}\subset E_{1,1}$\TODO
\end{exm}

\begin{prop}
\label{seed}
For $\ell,k\ge 0$, there is an isomorphism of fermionic groups $\varphi_{\ell,k}\colon E_{\ell+4,k}\xrightarrow\cong E_{\ell,k+4}$; restricted to $\Delta_{\ell+4,k}\subset E_{\ell+4,k}$, $\varphi_{\ell,k}$ defines an isomorphism $\Delta_{\ell+4,k}\cong\Delta_{\ell,k+4}$.
\end{prop}

Except for the bit about the diagonal subgroups, \cref{seed} is proven in~\cite[Corollary 2.31]{BottSpiral}, which constructs an explicit isomorphism $\varphi_{0,0}$ for the case $(\ell,k) = (0,0)$, then takes the fermionic tensor product with $E_{\ell,k}$ to define $\varphi_{\ell,k}$. The reader can check directly that this choice of $\varphi_{\ell,k}$ restricts to an isomorphism of the diagonal subgroups.

\begin{defn}
\label{EMon_defn}
    Let $\EMon$ be the commutative monoid formed from the set
    \begin{equation}
        \{E_{\ell,k}~|~\ell,k\geq 0\}/\cong
    \end{equation}
    equipped with the fermionic tensor product $\ftens$.
\end{defn}

\begin{lem}
\label{M_is_mod_4004}
The map $\varepsilon\colon \mathbb N^2\to\EMon$ defined by $(\ell,k)\mapsto E_{\ell,k}$ induces an isomorphism
    \begin{equation}
        \widetilde\varepsilon\colon \mathbb N^2/((4,0) - (0, 4)) \xrightarrow\cong \EMon.
    \end{equation}
\end{lem}

\begin{proof}
\Cref{seed} shows that $\epsilon$ factors through $\mathbb N^2/((4,0) - (0,4))$, so we have a surjective monoid homomorphism $\widetilde\varepsilon$ as in the lemma statement, and we are done if we can show $\widetilde\varepsilon$ is injective.

Let $\psi\colon \EMon\to \N\times\Z/8$ send $E_{\ell,k}\mapsto (\ell+k, \ell-k\bmod 8)$.
To see that this is a well-defined monoid homomorphism, we check that isomorphic fermionic groups $E_{\ell,k} \cong E_{\ell',k'}$ have equivalent images.
Isomorphic groups must have identical cardinality, so we have $\ell+k=\ell'+k'$.
% We will show that $\psi(E_{\ell,k})$ indeed detects the isomorphism type of $E_{\ell,k}$ as a fermionic group, so that it is a well-defined monoid homomorphism.
%That $\ell+k$ is an isomorphism invariant of $E_{\ell,k}$ is because it detects cardinality: $\# E_{\ell,k} = 2^{\ell+k+1}$.

For an additional invariant, we can use the following.
Given a fermionic group $G_f$, the \term{fermionic group algebra}~\cite[\S 7]{Ben88} $\R^f[G_f]$ is defined to be the quotient
\begin{equation}
    \R^f[G_f]\coloneqq \R[G_f]/(-1 + (-1)^F),
\end{equation}
where $(-1)^F$ is the unique nonzero element of the kernel of $G_f\to G_b$. The homomorphism $\theta\colon G_b\to\Z/2$ (see part~\eqref{thetapart} of \cref{fermgp}) gives $\R^f[G_f]$ the structure of a $\Z/2$-graded algebra over $\R$, and the Morita class as a $\Z/2$-graded algebra is an invariant of the isomorphism type of $G_f$.

Albuquerque--Majid~\cite{albuquerque2002clifford} show that $\R^f[E_{\ell,k}]$ is isomorphic as a $\Z/2$-graded algebra to the Clifford algebra $\Cl_{\ell,k}$.\footnote{Our Clifford algebra conventions match~\cite[Definition 1.47]{BottSpiral}: $\Cl_{\ell,k}$ is generated as an $\R$-algebra by $\ell+k$ odd, pairwise anticommuting elements of which $\ell$ square to $1$ and $k$ square to $-1$. In addition to Albuquerque--Majid~\cite{albuquerque2002clifford} , see also Thiang~\cite[Remark 6.4]{Thiang_2015}, Abłamowicz~\cite{ablamowicz2016cliffordalgebrasrelatedfinite}, Ichikawa--Tachikawa~\cite[Lemma 2.7]{IT23}, and~\cite[Corollary 1.67]{BottSpiral}.}
It is standard (see, e.g., \cite{atiyah_clifford_1963}) that $\Cl_{\ell,k}$ is $\Z/2$-graded Morita equivalent to $\Cl_{\ell',k'}$ if and only if $\ell-k\equiv\ell'-k'\bmod 8$. Thus $\ell-k\bmod 8$ is also an isomorphism invariant of $E_{\ell,k}$.

We want to show $\widetilde\varepsilon$ is injective, so we may as well show $\psi\circ\widetilde\varepsilon$ is injective. Suppose that $\ell,k,\ell',k'\in\N$ are such that $\ell+k = \ell'+k'$ and $\ell-k\equiv \ell'-k'\bmod 8$. Then there is some $x\in\Z$ such that $\ell-k = \ell'-k' + 8x$. Add this equation to $\ell+k = \ell'+k'$ to conclude that $2\ell=2\ell'+8x$ or $\ell = \ell'+4x$, and thus $k = k' - 4x$. That is, $(\ell,k)$ and $(\ell',k')$ can be identified by applying the relation $(4,0) = (0,4)$ $x$ times.
\end{proof}
%, if $E_{\ell,k}\cong E_{\ell',k'}$ then $\ell+k=\ell'+k'$.
%
%It thus remains to show that there are no other relations. 
%To do so, we will show that the pair $(\ell+k,\ell-k\bmod 8)\in\N\times\Z/8$ is an invariant of the isomorphism type of $E_{\ell,k}$ as a fermionic group. 
%This is equivalent to the fact that $((4,0) - (0,4))$ is the only relation in $\EMon$ because if one keeps $\ell+k$ fixed, the only way to change $\ell$ and $k$ without changing $\ell-k\bmod 8$ is to add $4$ to one of $\ell$ or $k$ and subtract $4$ from the other. Thus, if there were an additional relation in $\EMon$, one would obtain isomorphic fermionic groups with distinct values of $(\ell+k,\ell-k\bmod 8)$, which would be a problem.

%Therefore we are done if we can show that $E_{3,0}$, $E_{2,1}$, $E_{1,2}$, and $E_{0,3}$ are pairwise nonisomorphic fermionic groups.

\begin{prop}
For all $\ell,k,\ell',k'\ge 0$, the twisting data of $\Delta_{\ell+\ell',k+k'}$, as an element of the abelian group $[B\Z/2, B\O/B\Spin]$, is the sum of the twisting data of $\Delta_{\ell,k}$ and the twisting data of $\Delta_{\ell',k'}$. %\Nat{$\ell',k'$ not $\ell' + k'$?}
\end{prop}
\begin{proof}
Let $(\theta_1, \omega_1),(\theta_2,\omega_2)\in [BG_b, B\O/B\Spin]$ be the twisting data of two fermionic groups $G_f^{(1)},G_f^{(2)}$ whose bosonic groups are both isomorphic to $G_b$. The addition in the set of homotopy classes of maps to an $E_\infty$-space $A$ is defined as the pullback of the external sum (see \cref{slice_tensor}) by the diagonal map.

Unspooling this, the sum of $(\theta_1,\omega_1)$ and $(\theta_2,\omega_2)$ (again, we are using the non-product group law~\eqref{whit_twist}) is the pullback of the twisting data of $G_b\times G_b = (G_f^{(1)}\ftens G_f^{(2)})_b$ across the diagonal map $\Delta\colon G_b\to G_b\times G_b$.

Now let $G_b = \Z/2$, $G_f^{(1)} = \Delta_{\ell,k}$, and $G_f^{(2)} = \Delta_{\ell',k'}$. By definition, $\Delta_{\ell+\ell',k+k'}$ is the pullback
% https://q.uiver.app/#q=WzAsNCxbMCwwLCJcXERlbHRhX3tcXGVsbCtcXGVsbCcsaytrJ30iXSxbMCwxLCJcXFovMiJdLFsxLDEsIihcXFovMilee1xcZWxsK1xcZWxsJytrK2snfSJdLFsxLDAsIkVfe1xcZWxsK1xcZWxsJyxrK2snfSJdLFswLDNdLFswLDFdLFsxLDIsIlxcRGVsdGEiLDJdLFszLDJdLFswLDIsIiIsMSx7InN0eWxlIjp7Im5hbWUiOiJjb3JuZXIifX1dXQ==
\begin{subequations}
\begin{equation}
\begin{tikzcd}
	{\Delta_{\ell+\ell',k+k'}} & {E_{\ell+\ell',k+k'}} \\
	{\Z/2} & {(\Z/2)^{\ell+\ell'+k+k'}.}
	\arrow[from=1-1, to=1-2]
	\arrow[from=1-1, to=2-1]
	\arrow["\lrcorner"{anchor=center, pos=0.125}, draw=none, from=1-1, to=2-2]
	\arrow[from=1-2, to=2-2]
	\arrow["\Delta"', from=2-1, to=2-2]
\end{tikzcd}
\end{equation}
By construction, this map factors through $\Delta_{\ell,k}\ftens\Delta_{\ell',k'}$:
% https://q.uiver.app/#q=WzAsNixbMCwxLCJcXFovMiJdLFsxLDEsIihcXFovMilee1xcZWxsK1xcZWxsJytrK2snfSJdLFsxLDAsIlxcRGVsdGFfe1xcZWxsLGt9XFxmdGVuc1xcRGVsdGFfe1xcZWxsJyxrJ30iXSxbMCwwLCJcXERlbHRhX3tcXGVsbCtcXGVsbCcsaytrJ30iXSxbMiwwLCJFX3tcXGVsbCtcXGVsbCcsaytrJ30iXSxbMiwxLCIoXFxaLzIpXntcXGVsbCtcXGVsbCcraytrJ30iXSxbMCwxLCJcXERlbHRhIiwyXSxbMiwxXSxbMywwXSxbMywyXSxbMSw1LCIiLDAseyJsZXZlbCI6Miwic3R5bGUiOnsiaGVhZCI6eyJuYW1lIjoibm9uZSJ9fX1dLFsyLDRdLFs0LDVdLFszLDEsIiIsMSx7InN0eWxlIjp7Im5hbWUiOiJjb3JuZXIifX1dXQ==
\begin{equation}\label{triple_pullback}
\begin{tikzcd}
	{\Delta_{\ell+\ell',k+k'}} & {\Delta_{\ell,k}\ftens\Delta_{\ell',k'}} & {E_{\ell,k}\ftens E_{\ell',k'}} & {E_{\ell+\ell',k+k'}} \\
	{\Z/2} & {\Z/2\times\Z/2} & {(\Z/2)^{\ell+k}\times (\Z/2)^{\ell'+k'}} & {(\Z/2)^{\ell+\ell'+k+k'}}
	\arrow[from=1-1, to=1-2]
	\arrow[from=1-1, to=2-1]
	\arrow["\lrcorner"{anchor=center, pos=0.125}, draw=none, from=1-1, to=2-2]
    \arrow["\lrcorner"{anchor=center, pos=0.125}, draw=none, from=1-2, to=2-3]
	\arrow[from=1-2, to=1-3]
	\arrow[from=1-2, to=2-2]
	\arrow[from=1-3, to=2-3]
	\arrow["\Delta"', from=2-1, to=2-2]
	\arrow["{(\Delta, \Delta)}"', from=2-2, to=2-3]
    \arrow["\eqref{add_mon}"', equals, from=1-3, to=1-4]
    \arrow[equals, from=2-3, to=2-4]
    \arrow[from=1-4, to=2-4]
\end{tikzcd}\end{equation}
\end{subequations}
Apply the general discussion above to the leftmost square in~\eqref{triple_pullback} to see that, as we wanted to show, the twisting data of $\Delta_{\ell+\ell',k+k'}$ is the sum of the twisting data of $\Delta_{\ell,k}$ and $\Delta_{\ell',k'}$.
\end{proof}

\begin{cor}
\label{monoid_hom_cor}
The map sending $E_{\ell,k}$ to the twisting data of $\Delta_{\ell,k}$ is a homomorphism of commutative monoids $\delta\colon \mathbb M\to [B\Z/2, B\O/B\Spin]$.
\end{cor}
% [\TODO: define $\mathbb M$ earlier!]

\begin{cor}
\label{diag_mod_4}
If $\ell-k\equiv \ell'-k'\bmod 4$, there is an isomorphism of fermionic groups $\Delta_{\ell,k}\cong\Delta_{\ell',k'}$.
\end{cor}
\begin{proof}
%Since $H^1(B\Z/2; \Z/2)$ and $H^2(B\Z/2; \Z/2)$ are both isomorphic to $\Z/2$, $[B\Z/2, B\O/B\Spin]$ is an abelian group of order $4$, and therefore $4 = 0$ in it.
By \cref{seed}, the homomorphism $\mathbb M\to [B\Z/2, B\O/B\Spin]$ from \cref{monoid_hom_cor} identifies the images of $E_{4,0}$ and $E_{0,4}$. By \cref{elk_exm}, this homomorphism sends $E_{1,1}\mapsto 0$, hence the image of $E_{\ell,k}$ only depends on $\ell-k$. Combining these two facts, this homomorphism factors through the homomorphism $\mathbb M\to\Z/4$ sending $E_{\ell,k}\mapsto \ell-k\bmod 4$.
\end{proof}

\begin{exm}
\label{rank1_delta}
$E_{1,0}$ and $E_{0,1}$ are their own diagonal subgroups: by construction, $\Delta_{1,0}\to E_{1,0}$ is the pullback of the inclusion of the diagonal $\Z/2\subset (E_{1,0})_b \cong\Z/2$, and the pullback of an isomorphism is an isomorphism. Likewise for $E_{0,1}$. Therefore $\delta(E_{1,0})$ and $\delta(E_{0,1})$ are just the twisting data of these two fermionic groups, which we calculated in \cref{pin_exms}: $\delta(E_{1,0}) = (x, 0)$ and $\delta(E_{0,1}) = (x, x^2)$ (in the notation of \cref{{H1H2_twist}}).
\end{exm}

Since the isomorphism type of $\Delta_{\ell,k}$ only depends on $\ell-k\bmod 4$ by \cref{diag_mod_4}, we can write down the four classes.
\begin{enumerate}
    \setcounter{enumi}{-1}
    \item When $\ell-k\equiv 0\bmod 4$, $\Delta_{\ell,k}\cong\Delta_{1,1}\cong\Z/2\times\{\pm 1\}$, with trivial twisting data, as we saw in \cref{elk_exm}.
    \item When $\ell-k\equiv 1\bmod 4$, $\Delta_{\ell,k}\cong\Delta_{1,0}\cong\Pin_1^+$, as we saw in \cref{elk_exm}.
    \item When $\ell-k\equiv 2\bmod 4$, one can show that $\Delta_{\ell,k}\cong\Z/4$ with $\theta = 0$, e.g.\ by computing it inside $E_{2,0}$.
    \item When $\ell-k\equiv 3\bmod 4$, $\Delta_{\ell,k}\cong\Delta_{0,1}\cong\Pin_1^-$ as in \cref{elk_exm}.
\end{enumerate}

%\TODO: Define a symmetric monoidal $\infty$-category of fermionic groups (probably just groups $G$ with maps $BG\to B\O/B\Spin$), extract $\theta$ and $\omega$ and the fermionic tensor product. Introduce $E_{\ell,k}$ and its diagonal subgroup $\Delta_{\ell-k}$, construct the isomorphism $E_{4,0}\cong E_{0,4}$ and therefore introduce the commutative monoid
%\begin{equation}
%    \EMon \coloneqq \set{E_{\ell,k}\mid \ell,k\ge 0}/{\cong},
%\end{equation}
%with multiplication given by the fermionic tensor product. Thus, in particular, $\EMon\cong\mathbb N^2/((4,0) - (0, 4))$.

%%%%%%%%%%%%%%%%%%%%%
%%%%%%%%%%%%%%%%%%%%%

\section{The two periodicities}
\label{s:per}
In this section, we introduce the two $4$-periodicities that will be the main focus of this paper, as well as some background notions and definitions for them. These manifest as groups isomorphic to $\Z/4$, specifically the torsion subgroup of the $Q_1$-local stable Picard group of $\cA(1)$-modules (\cref{yu_thm}) and (a nice subgroup of) the twists of $\KO$-theory over $B\Z/2$ (\cref{second_Z4}), which we introduce in \S\ref{ss:first_Z4} and \S\ref{ss:second_Z4} respectively.

\subsection{Torsion in the Picard group of \texorpdfstring{$Q_1$}{Q1}-local \texorpdfstring{$\cA(1)$}{A(1)}-modules}
\label{ss:first_Z4}

Let $\cA$ denote the \term{mod $2$ Steenrod algebra}, the $\Z$-graded $\Z/2$-algebra of natural transformations $H^*(\bl;\Z/2)\to H^{*+k}(\bl;\Z/2)$ commuting with the suspension isomorphism. (The grading on $\cA$ is such that a natural transformation increasing the degree of a cohomology class by $k$ has degree $k$ in $\cA$.)
%the amount that such a natural transformation increases the degree by.) 
$\cA$ is generated by the \term{Steenrod squares} $\Sq^n$ in degree $n$ for each $n\ge 1$, subject to relations called \term{Adem relations}. For any space or spectrum $X$, $H^*(X;\Z/2)$ is naturally a $\Z$-graded $\cA$-module, essentially by definition of $\cA$.

The Cartan formula is the comultiplication map in a cocommutative Hopf algebra structure on $\cA$: see~\cite[Corollary 1 and Appendix 1]{milnor_steenrod_1958}. Thus the category of left $\cA$-modules acquires a symmetric monoidal structure with respect to tensor product over $\Z/2$.

%[\TODO: Hopf alg so there's a symm mon cat of $\cA$-modules] \Nat{does Hopf alg imply symmetric monoidal? I thought it only implies monoidal, unless the Hopf algebra is also cocommutative. See \cite[p. 187]{Mar83}. Milnor showed $\cA$ is a cocommutative Hopf algebra in Steenrod Algebra and its Dual, so $\text{Mod}_\cA$ is indeed symm mon, so we should cite that for the appropriate implications.}

For $n\ge 0$, let $\cA(n)$ denote the subalgebra of $\cA$ generated by $\Sq^k$ with $k\le 2^n$, which is a Hopf subalgebra. We will mostly focus on $\cA(1)$ in this note. Whenever we refer to $\cA$- or $\cA(n)$-modules, we implicitly mean in the graded sense: that a degree-$k$ element of $\cA$ or $\cA(n)$ acts on a homogeneous element $x$ in a way that raises the degree of $x$ by $k$.

\begin{defn}[{\cite{AM71, Mar83}}]
\label{marghom}
Let $Q_0\coloneqq\Sq^1$ and inductively define $Q_k\coloneqq [Q_{k-1}, \Sq^{2^k}]$. Thus $Q_k\in\cA(n)$ for $k\le n$. The Adem relations imply that $Q_k^2 = 0$, so for any $\cA(n)$-module $M$ and $k\le n$, one may regard $M$ as a chain complex with differential $Q_k$.

In this situation, the \term{$Q_k$-Margolis homology} of $M$, denoted $H_*(M; Q_k)$, is the homology of the above chain complex.
\end{defn}

Since $Q_k$ acts in a manner compatible with the grading, the $\Z$-grading of $M$ induces a $\Z$-grading on $H_*(M;Q_k)$. For any $k\le n$, $Q_k$-Margolis homology commutes with direct sums of $\cA(n)$-modules by definition, and satisfies a Künneth formula with respect to tensor products.
\begin{prop}[{Margolis~\cite[Chapter 19, Proposition
17(b)]{Mar83}}]
\label{margolis_kunneth}
For $\cA(n)$-modules $M$ and $N$, and $k\le n$, there is a natural isomorphism
\begin{equation}
\label{kunneth}
    H_*(M\otimes N; Q_k) \cong H_*(M; Q_k)\otimes H_*(N; Q_k).
\end{equation}
\end{prop}
Moreover, $H_*(\cA(n); Q_k) = 0$~\cite[\S 4]{AM71}.

\begin{defn}[{Adams--Margolis~\cite[Definition 4.1]{AM71}, Margolis~\cite[Chapter 23]{Mar83}}]
Let $n\in\N$, $k\le n$, and $f\colon M\to N$ be a map of bounded-below $\cA(n)$-modules.
\begin{enumerate}
    \item The map $f$ is a \term{stable isomorphism} if its kernel and cokernel are free $\cA(n)$-modules.
    \item The map $f$ is a \term{$Q_k$-local stable isomorphism} if is an isomorphism on $Q_k$-Margolis homology.
\end{enumerate}
\end{defn}
By~\cite[Chapter 14, Proposition 1]{Mar83}, $f\colon M\to N$ is a stable isomorphism if and only if there are free $\cA(n)$-modules $L_1$ and $L_2$ and an $\cA(n)$-module isomorphism
\begin{equation}
    \overline f\colon L_1\oplus M \xrightarrow\cong L_2\oplus N
\end{equation}
such that $\overline f|_M = (0, f)$.

\begin{prop}[{Wall~\cite[Theorem 2]{wall_characterization_1962}}]
The following are equivalent for a map $f\colon M\to N$ of bounded-below $\cA(1)$-modules.
\begin{enumerate}
    \item $f$ is a stable isomorphism.
    \item $f$ is both a $Q_0$-local and a $Q_1$-local stable isomorphism.
\end{enumerate}
\end{prop}

For $n>1$, Adams--Margolis~\cite[Theorem 4.2]{AM71} gave a similar characterization of stable isomorphisms, but using Margolis homology for a more general class of operations than $Q_k$. Nevertheless, since $H_*(\cA(n); Q_k) = 0$ for any $n\in\N$ and $k\le n$, every stable isomorphism is a $Q_k$-local stable isomorphism.

For any $k\le n$, it is immediate that the $Q_k$-local equivalences contain all $\cA(n)$-module isomorphisms and satisfy a two-out-of-three property. Following \cite[Proposition 4.1.7.4]{HA}, to obtain a good theory of localizations, we also need that $Q_k$-local equivalences are stable with respect to tensor products, which follows from the Künneth theorem (\cref{margolis_kunneth}). Thus we conclude:

\begin{thm}[{Ricka~\cite[Proposition 2.6]{ricka_local_2016}}]
\label{local_stmod}
There is a symmetric monoidal stable $\infty$-category $\StMod_{\cA(n)}[Q_k^{-1}]$ which is the localization of the category of bounded-below $\cA(n)$-modules at the $Q_k$-local equivalences, and the localization functor $L_k\colon \Mod_{\cA(n)}\to\StMod_{\cA(n)}[Q_k^{-1}]$ is symmetric monoidal.
\end{thm}

It is also possible to construct the \term{stable module category} $\StMod_{\cA(n)}$, the localization of the category of $\cA(n)$-modules at the stable isomorphisms, as a symmetric monoidal $\infty$-category in a similar way. See Ricka~\cite[Proposition 1.9]{ricka_local_2016} and Yang~\cite[\S 3.5]{yang_thesis}. We will not need this so will not say much more about it.
%\begin{defn}
%Let $n\in\N$ and $k\le n$.
%\begin{enumerate}
%    \item The \term{stable $\cA(n)$-module category}, $\StMod_{\cA(n)}$, is the symmetric monoidal $\infty$-category obtained by localizing at the stable isomorphisms.
%    \item The \term{$Q_k$-local stable $\cA(n)$-module category}, denoted , is defined analogously for the $Q_k$-local stable isomorphisms.
%\end{enumerate}
%\end{defn}%
%
%In both cases the symmetric monoidal structure is induced from the tensor product of $\cA(n)$-modules, itself coming from the cocommutative Hopf algebra structure on $\cA(n)$ as described above. 
\begin{defn}
Let $(\cat C, \otimes, \boldsymbol 1)$ be a symmetric monoidal category.
\begin{enumerate}
    \item An object $x\in\cat C$ is \term{invertible} if there is another object $x^{-1}\in\cat C$ and an isomorphism $x\otimes x^{-1}\cong \boldsymbol 1$.
    \item The \term{Picard group} of $\cat C$, denoted $\Pic(\cat C)$, is the abelian group of isomorphism classes of invertible objects, with unit the isomorphism class of $\boldsymbol 1$ and the group operation given by $\otimes$.
\end{enumerate}
If $\cat C$ is a symmetric monoidal $\infty$-category, we define $\Pic(\cat C)$ to be the Picard group of the homotopy category of $\cat C$.
\end{defn}
See~\cite{hopkins_constructions_1994, may_picard_2001, van_der_meer_endotrivial_2021} for more background on Picard groups, and some example computations.
\begin{exm}
The Picard groups of $\StMod_{\cA(n)}$ are known for all $n$.
\begin{enumerate}
    \item Adams--Priddy~\cite[Theorem 3.7]{AP76} constructed an isomorphism $\Pic(\StMod_{\cA(1)})\cong\Z^2\oplus\Z/2$; the unique nonzero torsion class is that of the \term{Joker} $\Sigma^{-2}(\cA(1)/\Sq^3)$.
    \item Bhattacharya--Ricka~\cite[Theorem 4.6]{bhattacharya_stable_2017} showed that $\Pic(\StMod_{\cA(2)})\cong\Z^2$.
    \item Pan--Yan~\cite[Corollary 4.8]{pan_stable_2026} showed that for all $n\ge 3$, $\Pic(\StMod_{\cA(n)})\cong\Z^2$.
    \item Gheorghe--Isaksen--Ricka~\cite[Theorem 5.6]{gheorghe_picard_2018} and Gao--Li~\cite[Theorem 4.3]{gao_picard_2024} prove that the Picard groups of the $\C$-, resp.\ $\R$-motivic analogues of $\StMod_{\cA(1)}$ are both isomorphic to $\Z^4$.
    \qedhere
\end{enumerate}
\end{exm}
We will be interested in the $Q_1$-local Picard group of $\cA(1)$.
\begin{thm}[{Yu~\cite[Theorem 2.1]{Yu95}}]
\label{yu_thm}
There is an isomorphism \begin{equation}
    y\colon \Pic(\StMod_{\cA(1)}[Q_1^{-1}])\xrightarrow\cong \Z\oplus\Z/4
\end{equation}
such that:
\begin{itemize}
    \item $y^{-1}((1, 0))$ is the $Q_1$-localization of $\Sigma\Z/2$, and
    \item $y^{-1}((0, 1))$ is $T_1\coloneqq \Sigma^{-1} H^*((B\Z/2)^{\sigma-1};\Z/2)$ (which is $Q_1$-local).
\end{itemize}
The isomorphism $y$ induces an isomorphism from the torsion summand of $\Pic(\StMod_{\cA(1)}[Q_1^{-1}])$ to $\Z/4$; the module $T_i$ in \cref{Q1_torsion} is invertible in $\StMod_{\cA(1)}[Q_1^{-1}]$ and $y(T_i) = i\in\Z/4$.
%is sent to $i\in\Z/4$ under that isomorphism.
%under that isomorphism, the module $T_i$ in~\cref{Q1_torsion} represents $i
\end{thm}
See also Bruner~\cite[Theorem 7.1 and Proposition 8.2]{Bru14}, who gave another proof of \cref{yu_thm} and also showed $\Pic(\StMod_{\cA(1)}[Q_0^{-1}])\cong\Z$. As far as we know, these are the only values of $k$ and $n$ for which the isomorphism type of $\Pic(\StMod_{\cA(n)}[Q_k^{-1}])$ is known.
In \cref{Q1_torsion}, we draw $Q_1$-local $\cA(1)$-modules $T_1$, $T_2$, $T_3$, and $T_4$ representing the four torsion classes in $\Pic(\StMod_{\cA(1)}[Q_1^{-1}])$.

The first of our two $\Z/4$s is the torsion subgroup of $\Pic(\StMod_{\cA(1)}[Q_1^{-1}])$. It can be interpreted as a periodicity in the following way: if $M$ is a bounded-below, $Q_1$-local $\cA(1)$-module, then the isomorphism type of $M\otimes T_1^{\otimes \ell}$, modulo free $\cA(1)$-module summands, is $4$-periodic in $\ell$. This fact was used to make computations of the $\ko$-homology of $(B\Z/2)^n$ in~\cite[Theorem A]{Yu95}, and of certain twisted $\ko$-homology groups over $(B\Z/2)^n$ in~\cite[\S 5]{BottSpiral}.

A corollary of \cref{yu_thm} is that there is a preferred projection from $\Pic(\StMod_{\cA(1)}[Q_1^{-1}])$ onto its torsion subgroup.
\begin{lem}
\label{desuspend}
The function
\begin{equation}
\varsigma\colon \Pic(\StMod_{\cA(1)}[Q_1^{-1}]) \longrightarrow
\Pic(\StMod_{\cA(1)}[Q_1^{-1}])_{\mathrm{tors}}
\end{equation}
defined by suspending or desuspending a $Q_1$-locally invertible $\cA(1)$-module until it is stably isomorphic to one of the modules $T_i$ is a group homomorphism.
\end{lem}
\begin{proof}
After applying the isomorphisms
\begin{equation}
\begin{gathered}
y\colon \Pic(\StMod_{\cA(1)}[Q_1^{-1}])\xrightarrow\cong\Z\oplus\Z/4\\y|_{\mathrm{tors}}\colon \Pic(\StMod_{\cA(1)}[Q_1^{-1}])_{\mathrm{tors}}\xrightarrow\cong\Z/4
\end{gathered}
\end{equation}
from \cref{yu_thm}, $\varsigma$ is identified with the function $(a, b)\mapsto b$, which is a homomorphism.
\end{proof}

\begin{figure}[h!]
\begin{tikzpicture}[scale=0.6]
\begin{scope}[every node/.style = {font=\tiny}]
    \draw[gray!50!white] (-5.5, 0) -- (12, 0);

		\clip (-6.25, -3.25) rectangle (11.5, 9.5);
% N0
        \foreach \y in {-2, 0, ..., 8} {
            \node[gray] at (-6, \y) {$\y$};
        }
        \foreach \y in {0, 2,..., 10} {
		%	\node at (-1.5, \y+1) {$\y$};
			\sqone(-4.25, \y-1);
		}
		\foreach \y in {0, ..., 10} {
			\tikzpt{-4.25}{\y-1}{}{};
		}
		\foreach \y in {0, 4, 8} {
			\sqtwoR(-4.25, \y-1);
		}
		\foreach \y in {3, 7} {
			\sqtwoL(-4.25, \y-1);
		}

% N1
		\foreach \y in {-1, 1,..., 9} {
			%\node at (-1.5, \y) {$\y$};
			\sqone(0, \y);
		}
		\foreach \y in {-1, ..., 9} {
			\tikzpt{0}{\y}{}{};
		}
		\foreach \y in {0, 4, 8} {
			\sqtwoL(0, \y);
		}
		\foreach \y in {1, 5, 9} {
			\sqtwoR(0, \y);
		}

% N2
    \foreach \y in {0, 1, 3, 4} {
		  \tikzpt{4}{\y-2}{}{};
	}
	\foreach \y in {1, ..., 11} {
		\tikzpt{5.5}{\y-2}{}{};
	}
	\sqone(4, -2);\sqone(4, 1);
	\sqtwoL(4, -1);
	\sqtwoCR(4, -2);
	\sqtwoCL(5.5, 0);
	\foreach \y in {1, 3, ..., 11} {
		\sqone(5.5, \y-2);
	}
	\begin{scope}
		\foreach \y in {1, 5, 9} {
	%		\clip (0, 0) rectangle (2.5, 10.75);
			\sqtwoR(5.5, \y-2);
		}
	\end{scope}
	\foreach \y in {4, 8} {
		\sqtwoL(5.5, \y-2);
	}

% N3
		\foreach \y in {0, 2,..., 12} {
			\sqone(9.5, \y-3);
		}
		\foreach \y in {0, ..., 12} {
			\tikzpt{9.5}{\y-3}{}{};
		}
        \sqtwoL(9.5, -3);
		\foreach \y in {4, 8, 12} {
			\sqtwoR(9.5, \y-3);
		}
		\foreach \y in {3, 7, 11} {
			\sqtwoL(9.5, \y-3);
		}
        \tikzpt{11}{0}{}{};
        \tikzpt{11}{1}{}{};
        \sqone(11, 0);
        \sqtwoCR(9.5, -2);
        \sqtwoCR(9.5, -1);
    \end{scope}
        \node at (-4.25, -4) {$T_0$};
        \node at (0, -4) {$T_1$};
        \node at (4.75, -4) {$T_2$};
        \node at (9.5, -4) {$T_3$};
\end{tikzpicture}
\caption{By \cref{yu_thm}, the torsion subgroup of $\Pic(\protect\StMod_{\cA(1)}[Q_1^{-1}])$ is isomorphic to $\Z/4$. In this figure, the $\cA(1)$-module $T_i$ is in the class $i\in\Z/4\cong \Pic(\protect\StMod_{\cA(1)}[Q_1^{-1}])_{\mathrm{tors}}$. Each module $T_i$ is the unique $\cA(1)$-module with this property, up to isomorphism, which is $Q_1$-local and has no free $\cA(1)$-module summands. Though we have but a finite amount of space, each $T_i$ continues upward in the same $4$-periodic pattern.}
\label{Q1_torsion}
\end{figure}

%\TODO:
%Introduce the $Q_i$-local stable module category of $\cA(n)$ and Yu's work~\cite{Yu95} calculating the Picard group for $i = 1$, $n = 1$. Discuss the torsion summand, giving a $\Z/4$ periodicity

\subsection{Smith periodicity}
\label{ss:second_Z4}
In this subsection, we discuss our second $\Z/4$, the group $[B\Z/2, B\O/B\Spin]$ (\cref{second_Z4}). We discuss how this applies to Picard-group-graded twisted $\KO$-theory of $B\Z/2$ as well as a relation to families of Smith homomorphisms in bordism theory.
%\begin{itemize}
%    \item $B\O/B\Spin$ gives you geometrically meaningful twists of spin bordism (hence also twisted $\ko$) -- define twisted tangential strs more generally AND/OR give motivation for spin structures
    %\item Introduce the periodicity coming from the $E_\infty$-%structure on $B\O/B\Spin$, and why it is interesting (e.g.\ $\Z\times\Z/4$-graded twisted $\ko$-cohomology)
    %\item We should probably prove that $[B\Z/2, B\O/B\Spin]\cong \Z/4$ (since $B\O/B\Spin = K(\Z/2, 1)\times K(\Z/2, 2)$ as spaces, there are $4$ homotopy classes of maps. And in this group, we can show $2\ne 0$ trivially)
%\end{itemize}
\begin{notn}
\label{T_defn}
Recall the $E_\infty$-space $B\O/B\Spin$ from \cref{bo_bspin}, and for an $E_\infty$-ring spectrum $R$,\footnote{The construction works more generally for $R$ an $E_1$-ring spectrum: see Ando--Blumberg--Gepner--Hopkins--Rezk~\cite{ABGHR14a, ABGHR14b}. We stick to only the level of generality we need.} let $\GL_1(R)$ be the $\infty$-group of units of $R$, as defined in May~\cite[\S III.2]{MQRT77}.
May--Quinn--Ray~\cite[Lemma IV.2.6]{MQRT77} construct an $E_\infty$-map $T\colon B\O/B\Spin\to B\GL_1(\MTSpin)$, where $\MTSpin$ is the Thom spectrum whose homotopy groups are the spin bordism groups. 
Using the Atiyah--Bott--Shapiro maps $\MTSpin\to\ko\to\KO$~\cite{atiyah_clifford_1963}, which Joachim~\cite{Joa04} showed admit $E_\infty$-structures, we then get $E_\infty$-maps $T\colon B\O/B\Spin\to B\GL_1(\ko)\to B\GL_1(\KO)$.\footnote{By definition, the space $B\GL_1(R)$ only depends on $\Omega^\infty(R)$, so the map $B\GL_1(\ko)\to B\GL_1(\KO)$ is in fact an equivalence of $E_\infty$-spaces.}
\end{notn}

\begin{prop}
\label{gl1ko_split}
$T\colon B\O/B\Spin\to B\GL_1(\ko)$ admits a section, and there is an equivalence of $E_\infty$-spaces
\begin{equation}
    B\GL_1(\ko)\simeq B\O/B\Spin \times B(B\SO_\otimes),
\end{equation}
where $B\SO_\otimes$ is the space $B\SO$ with $E_\infty$-structure induced from the tensor product of vector bundles.
\end{prop}

Madsen--Snaith--Tornehave~\cite[Lemma 1.13]{madsen_infinite_1977} prove $B\SL_1(\ko)\simeq B\SO/B\Spin\times B(B\SO_\otimes)$, and the proof of \cref{gl1ko_split} is completely analogous---in particular, \cref{gl1ko_split} is not a new result. See also Antieau--Gepner--Gómez~\cite[\S 1]{antieau_actions_2014}, who prove a version of \cref{gl1ko_split} but do not discuss $E_\infty$-structures.

We care about $B\GL_1(R)$ because it grades twisted $R$-cohomology, as follows.
\begin{defn}[{Ando--Blumberg--Gepner~\cite[Definition 1.2]{ABG18}}]
Let $R$ be an $E_\infty$-ring spectrum and $f\colon X\to B\GL_1(R)$ be a map. The \term{$f$-twisted $R$-cohomology groups} of $X$ are
\begin{equation}
    R^{*-f}(X) \coloneqq \pi_{-*}\mathrm{Map}_R(M^R f, X),
\end{equation}
where $M^R f$ denotes the $R$-module Thom spectrum of $f$, defined in~\cite[Definition 2.20]{ABGHR14a}.
\end{defn}

%[\TODO: return to this and the ``$*+f$'' part once we've gotten further in (i.e. check the twist sign convention against $K^0(X^V) = K^{-V}(X)$)]

Donovan--Karoubi~\cite{donovan_graded_1970} defined twisted $\KO$-groups associated to a space $X$, $a\in H^1(X;\Z/2)$, and $b\in H^2(X;\Z/2)$. See also \cite{karoubi_algebres_1968}.

\begin{prop}[{Hebestreit--Sagave~\cite{HJ20}}]
Donovan--Karoubi's twisted $\KO$-theory groups are naturally isomorphic to the twisted $\KO$-cohomology groups, in the Ando--Blumberg--Gepner sense, of the map
\begin{equation}
    X \xrightarrow{(a,b)} K(\Z/2, 1)\times K(\Z/2, 2) \xrightarrow[\eqref{bobsp_split}]{\phi^{-1}} B\O/B\Spin \xrightarrow{T} B\GL_1(\KO).
\end{equation}
\end{prop}

Since $\KO$ is an $E_\infty$-ring spectrum, the $\KO$-cohomology groups of a space $X$ assemble into a $\Z$-graded ring, and in fact a $\Z/8$-graded ring thanks to Bott periodicity. For a fixed twist $(a,b)$, the twisted $\KO$-groups of $X$ do not in general admit a graded ring structure, only a $\KO^*(X)$-module structure.\footnote{However, in $K$- and $\KO$-\emph{homology}, there can be ring structures induced by the Pontrjagin product: see Douglas~\cite[\S 2.2]{douglas_twisted_2006} and Freed--Hopkins--Teleman~\cite[\S 3.6]{freed_loop_2011}.} This is similar to the structure present (or absent) in cohomology with local coefficients.

On the other hand, if one considers all twists at once, it is possible to multiply: given maps $f_1,f_2\colon X\to B\GL_1(\KO)$ and classes $x_i\in \KO^{n_i - f_i}(X)$, there is a product $x_1x_2\in \KO^{n_1 + n_2-f_1-f_2}(X)$, such that~\cite{donovan_graded_1970}
\begin{equation}
    \KO^\bullet(X)\coloneqq\bigoplus_{f\in [X, B\GL_1(\KO)]}  \KO^{*-f}(X)
\end{equation}
has the data of a $\Z/8\times [X, B\GL_1(\KO)]$-graded ring. The subring in gradings $\Z/8\times\set 0$ is the untwisted $\KO$-cohomology of $X$ with its usual ring structure. This is a cousin of the $\Z\times H^1(X;\Z/2)$-grading on ordinary cohomology with local coefficients: see~\cite{cadek_cohomology_1999, lastovecki_cohomology_2005, wendt_chow_2024, brazelton_chow_2025, DYY23} for computations using this graded ring structure.

By \cref{gl1ko_split}, the $[X, B\GL_1(\KO)]$-grading splits into two pieces: a grading by $[X, B\O/B\Spin]$ and one by $[X, B(B\SO_\otimes)]$. The latter is in general quite mysterious,\footnote{See Dadarlat--Pennig~\cite{dadarlat_unit_2015} and Pennig~\cite{pennig_noncommutative_2016} for geometric models of the analogue of $[X, B(B\SO_\otimes)]$ for complex $K$-theory.} so we will focus on the former.
%[\TODO: do the higher twists actually vanish for $B\Z/2$? if so, say sth like ``for our application they don't matter'']

\begin{lem}
\label{second_Z4}
$[B\Z/2, B\O/B\Spin]\cong\Z/4$.
\end{lem}
Thus the twisted $\KO$-cohomology groups of $B\Z/2$ form a $\Z/8\times\Z/4$-graded ring. It would be interesting to compute this graded ring.
\Cref{second_Z4} was stated without proof in~\cite[\S 1.2.3]{DY23} and~\cite[Prop. C.37 (1)]{lin_topological_2024}.
\begin{proof}
By \cref{bobsp_split}, there is a homotopy equivalence $B\O/B\Spin\simeq K(\Z/2, 1)\times K(\Z/2, 2)$, so there is a bijection \emph{that is not necessarily a group isomorphism}
\begin{equation}
    [B\Z/2, B\O/B\Spin]\approx [B\Z/2, K(\Z/2, 1)\times K(\Z/2, 2)] = H^1(B\Z/2;\Z/2)\times H^2(B\Z/2; \Z/2).
\end{equation}
Thus $[B\Z/2, B\O/B\Spin]$ is a group with four elements, so to finish the proof, it suffices to exhibit an element $x\in [B\Z/2, B\O/B\Spin]$ with $2x\ne 0$. For this, we use the fact that, since the map $q\colon B\O\to B\O/B\Spin$ is $E_\infty$ and the composition
\begin{equation}
    B\O\to B\O/B\Spin \xrightarrow[\eqref{bobsp_split}]{\phi} K(\Z/2, 1)\times K(\Z/2, 2)
\end{equation}
sends $V\mapsto (w_1(V), w_2(V))$ (see~\cite[(1.42)]{DY23}, and the related~\cite[Lemma 7]{donovan_graded_1970}), we can add classes in $[B\Z/2, B\O/B\Spin]$ by passing through $\phi$ and using the Whitney sum formula: given $(a_1, b_1)$ and $(a_2,b_2)$ in $H^1(X;\Z/2)\times H^2(X;\Z/2)$, their sum as elements of $[X, B\O/B\Spin]$ is
\begin{equation}
\label{whitney_not}
    (a_1, b_1) \oplus (a_2, b_2) = (a_1 + a_2, b_1 + a_1a_2 + b_2).
\end{equation}
Thus if $x\in H^1(B\Z/2; \Z/2)$ is the unique nonzero class, the sum of the twists $(x,0)$ and $(x,0)$ is $(0, x^2)$, which is nonzero.
\end{proof}
\begin{rem}
One can show that $[B\Z/2, B(B\SO_\otimes)]\cong\Z/2$. Thus we can consider a non-null-homotopic map $f\colon B\Z/2\to B(B\SO_\otimes)$ and then map to $B\GL_1(\ko)$. Such a map is an example of a ``higher twist'' of $\ko$, and it would be interesting to better understand how it interacts with the $\Z/4$ we discuss in this section, either from the point of view of twisted $\KO$-theory or from that of twisted spin structures. Brooks~\cite{brook_higher_2020} makes some related computations of higher twists of complex $K$-theory.

One can compute $[B\Z/2, B(B\SO_\otimes)]$ in a few ways. Since $B\Z/2$ is $2$-local, we may $2$-localize $B(B\SO_\otimes)$; then, Adams--Priddy~\cite[Corollary 1.4]{AP76} identify $[\bl, B(B\SO_\otimes)_{(2)}]$ with $\tau_{\ge 2}\ko_{(2)}$, where $\tau_{\ge 2}$ denotes the Postnikov $1$-connected cover. So we want to understand $(\tau_{\ge 2}\ko)^1(B\Z/2)$; here are two ways forward.
\begin{itemize}
    \item The cofiber sequence $\tau_{\ge 2}\ko\to \ko\to\tau_{\le 1}\ko$ induces a long exact sequence in the cohomology of $B\Z/2$. Bruner--Greenlees~\cite[\S A.9]{BG10} calculate $\ko^*(B\Z/2)$, and one can attack $\tau_{\le 1}\ko^*(B\Z/2)$ with the Atiyah--Hirzebruch spectral sequence.
    \item Alternatively, since $\tau_{\ge 2}\ko$ is a $\ko$-module, there is a natural isomorphism
    \begin{equation}\label{Adams_ko_coh}
        \tau_{\ge 2}\ko^n(X) \coloneqq \pi_{-n}\mathrm{Map}(X, \tau_{\ge 2}\ko) \cong\pi_{-n}\mathrm{Map}_\ko(\ko\wedge X, \tau_{\ge 2}\ko).
    \end{equation}
    Having rephrased this question internal to the category of $\ko$-module spectra, we can finish it with Baker--Lazarev's $\ko$-module Adams spectral sequence~\cite{BL01}.\footnote{Boardmann~\cite[Theorem 15.6]{Boa99} showed that half-plane Adams spectral
sequences such as~\eqref{Adams_ko_coh} are conditionally convergent, since $\ko$ and $B\Z/2$ are finite type, and thus one must in general be careful with convergence. However, this particular spectral sequence collapses, so it strongly converges (\textit{ibid.}, \S 7) and everything is OK.} See~\cite{Bak20, DY23, kuroda_computations_2025, BottSpiral, DYY} for some examples of computations with this spectral sequence over $\ko$; in particular, using standard techniques in $\cA(1)$-modules one can show the spectral sequence computing $[B\Z/2, B(B\SO_\otimes)]$ collapses to imply the claim.
\end{itemize}
By contrast, Gómez~\cite{gomez_nonexistence_2010} showed that if $G$ is \emph{any} compact Lie group, there are no higher twists of complex $K$-theory over $BG$.
\end{rem}
The above story using maps to $B\GL_1(R)$ to define a graded ring structure on twisted $R$-cohomology is closely related to the extended gradings on (equivariant) parametrized cohomology discussed in work of Costenoble, Hudson, Tilson, and Waner~\cite{costenoble_equivariant_1992, costenoble_equivariant_2016, costenoble_ordinary_2017, costenoble_equivariant_2022, costenoble_geometric_2023, costenoble_graded_2023, costenoble_algebraic_2024, costenoble_equivariant_2024a, costenoble_equivariant_2025a, costenoble_equivariant_2025b, costenoble_equivariant_2025, costenoble_equivariant_2026} and Beaudry--Lewis--May--Petersen--Tatum~\cite{beaudry_guide_2024,beaudry_parametrized_2026}, as well as to the \term{permutation twisted cohomology} of a finite group as introduced by Balmer--Gallauer~\cite[\S 10]{balmer_geometry_2025} and also studied by Gómez and Miller~\cite{miller_permutation_2025, gomez_classification_2026, gomez_permutation_2026}.

We finish this subsection with a more geometric description of the isomorphism in \cref{second_Z4}.
Recall that a vector bundle admits a spin structure if its first and second Stiefel--Whitney classes vanish.
Twists of spin structure quantify exactly how a vector bundle may \emph{fail} to be spin.

The following definition, in the special case $a = 0$, is due to B.L. Wang~\cite[Definition 8.2]{Wan08}. It is also closely related to Kreck's definition of a normal $1$-smoothing~\cite[\S 2]{kreck_surgery_1999}.

\begin{defn}
Let $X$ be a space, $a\in H^1(X;\Z/2)$, and $b\in H^2(X;\Z/2)$. An \term{$(X,a, b)$-twisted spin structure} on a vector bundle $E\to M$ is the data of a map $f\colon M\to X$ and trivializations of $w_1(E) + f^*(a)$ and $w_2(E) + f^*(a)^2 + f^*(b)$.
\end{defn}

By a \term{trivialization} of a cohomology class, we mean a null-homotopy of the representing map to an Eilenberg--Mac Lane space. Also, for $(X, a, b)$ fixed, the $(X,a,b)$-twisted spin structures on a vector bundle $E\to M$ define a homotopy type; we consider two $(X,a,b)$-twisted spin structures to be equivalent if they belong to the same path component.

Using \cref{bobsp_split}, we will identify $(a,b)$ with a homotopy class of maps $X\to B\O/B\Spin$.

\begin{lem}
\label{twist_descent}
Let $V\to X$ be a vector bundle, $a\in H^1(X; \Z/2)$, and $b\in H^2(X;\Z/2)$. Let $M$ be a closed manifold whose tangent bundle has an $(X, a, b)$-twisted spin structure induced by a map $f\colon M\to X$. If $N\subset M$ is a closed submanifold whose mod $2$ homology class is Poincaré dual to the mod $2$ Euler class of $f^*(V)$, then $N$ has a canonical twisted spin structure for the sum
\begin{equation}
    (a, b)\oplus q(V - \mathrm{rank}(V)) \in [X, B\O/B\Spin].
\end{equation}
\end{lem}

See also~\cite[Lemma A.2]{DYY23}.
Here $q\colon B\O\to B\O/B\Spin$ is the quotient map. Note that we subtracted off the rank of $V$ since $B\O$ classifies rank-zero virtual vector bundles.

\begin{proof}
Since the fundamental class of $N$ is Poincaré dual to $e(V)$, there is an isomorphism $TM|_N\cong TN \oplus f^*V|_N$. Then use the Whitney sum formula.
\end{proof}

\Cref{twist_descent} organizes different kinds of twisted spin structure over the same space $X$ into families related by iteratively passing from $M$ to $N$.

\begin{exm}
\label{explicit_twists}
\Cref{second_Z4} implies that for twists of spin structure over $B\Z/2$, we obtain a $4$-periodic family. Specifically, if $\sigma\to B\Z/2$ denotes the tautological bundle, we obtain the following four structures. Under the isomorphism $[B\Z/2, B\O/B\Spin]\xrightarrow\cong\Z/4$ from \cref{second_Z4}, entry $(i)$ in the list below maps to $i\in\Z/4$.
\begin{enumerate}
    \setcounter{enumi}{-1}
    \item We begin with the trivial twist, where $a = 0$ and $b = 0$. A $(B\Z/2, 0, 0)$-twisted spin structure on $E\to M$ is the data of a spin structure on $E$ and a principal $\Z/2$-bundle $P\to M$. This is sometimes called a spin$\times\Z/2$ structure.
    
    \item The next entry in the family is a $(B\Z/2, x, 0)$-twisted spin structure, where $x\coloneqq w_1(\sigma)$. By Peterson~\cite[\S 7]{Pet68}, this is equivalent to a \pinm structure.
    
    \item Applying~\eqref{whitney_not}, our next twist is $(B\Z/2, 0, x^2)$. By Campbell~\cite[\S 7.9]{Cam17}, this kind of twisted spin structure is a $G$-structure for $G = \Spin\times_{\set{\pm 1}}\Z/4$.
    
    \item The next twist is $(B\Z/2, x, x^2)$-twisted spin structure, which by Stolz~\cite[\S 8]{Sto88} is a \pinp structure.
\end{enumerate}
If one applies~\eqref{whitney_not} once more, one gets back to the trivial twist, illustrating that this family is $4$-periodic.
See \cite[Example 7.8]{debray_smith_2024} and references therein.
In particular, this family admits a physical interpretation in terms of dimensional reductions of quantum field theories; see e.g.\ \cite{kapustin_fermionic_2015,TY19,HKT19,COSY,WWZ19,debray_long_2024,debray_smith_2024}.
\end{exm}

\begin{rem}
The periodicity of twists induced by a cyclic subgroup of $[X, B\O/B\Spin]$ lifts to a periodicity of iterated \term{Smith homomorphisms}, which are roughly speaking the effects on bordism classes of the construction described in \cref{twist_descent}. This is the context in which this fourfold periodicity for $X = B\Z/2$ is studied in the references cited above. See~\cite[\S 6]{debray_smith_2024} for more information.

% The passage from the bordism class of $M$ to the class of hte submanifold $N$ described in \cref{twist_descent} defines a \textit{Smith homomorphism} of bordism groups. 
Although we discuss them here in terms of the dual to the \textit{cohomology} Euler class, these homomorphisms are more naturally understood as representing the Euler class in \textit{cobordism}; see e.g. ~\cite[Appendix B]{debray_smith_2024}.
\end{rem}

%Review the story from~\cite{debray_smith_2024} (and other refs before) on Smith families and the $\Z/4$ periodicity of twists of spin bordism or $\ko$ over $B\Z/2$

%%%%%%%%%%%%%%%%%%%%%
%%%%%%%%%%%%%%%%%%%%%

\section{From fermionic groups to \texorpdfstring{$Q_1$}{Q1}-local \texorpdfstring{$\cA(1)$}{A(1)}-modules}
\label{s:assemble}
To prove \cref{Q1mainthm}, we need to pass from fermionic groups to $\Pic(\StMod_{\cA(1)}[Q_1^{-1}])$. We will do this by the following chain of symmetric monoidal functors, whose notation we will explain and/or recall below:
\begin{equation}
\label{many_functors}
    \Phi\colon \cat{FermGrp} \xrightarrow{\mathbf T}
    \cat{Top}^{\mathit{ft}}/_{B\GL_1(\ko)} \xrightarrow M
    %\cat{Mod}_{\MTSpin} %\xrightarrow{\widehat A} 
    \cat{Mod}_{\ko}^{\mathit{ft}} \xrightarrow{H_\ko^*} \cat{Mod}_{\cA(1)} ^{\mathit{ft}}\xrightarrow[\eqref{local_stmod}]{L_1} \StMod_{\cA(1)}[Q_1^{-1}].
\end{equation}
\begin{enumerate}
    \item Recall from \cref{FGdefn} that $\cat{FermGrp}$ is a full subcategory of $\cat{Top}/_{B\O/B\Spin}$. The functor $\mathbf T$ composes the map $BG\to B\O/B\Spin$ with the $E_\infty$-map $T\colon B\O/B\Spin\to B\GL_1(\ko)$ from \cref{T_defn}. By \cref{slice_tensor}, the $E_\infty$-structure on $T$ induces a symmetric monoidal structure on $\mathbf T$.

    The superscript $\mathit{ft}$ means finite type: since we have restricted to finite fermionic groups, their classifying spaces are finite type.
    
    \item The functor $M$ is Ando--Blumberg--Gepner--Hopkins--Rezk's Thom spectrum functor~\cite[Definition 2.20]{ABGHR14a}, which by Ando--Blumberg--Gepner~\cite[Theorem 1.6]{ABG18} is symmetric monoidal. Moreover, the Thom isomorphism theorem (specifically its incarnation in~\cite[Lemma 2.24]{DY23}) can be used to show that $M$ sends finite-type spaces over $B\GL_1(\ko)$ to finite-type $\ko$-modules.
    
    \item For a $\ko$-module $M$, let $H_\ko^*(M)\coloneqq \pi_{-*}\mathrm{Map}_\ko(M, H\Z/2)$~\cite[Definition IV.1.7]{EKMM97}, where $H\Z/2$ carries the $E_\infty$-$\ko$-algebra structure given by Postnikov truncation $\ko\to H\Z$ followed by reduction mod $2$. This is naturally a module over $H_\ko^*(H\Z/2)$, which Baker~\cite[Theorem 5.1]{Bak20} shows has a unique isomorphism to $\cA(1)$. Like for ordinary cohomology, $H_\ko^*$ of a finite type $\ko$-module is a finite-type $\cA(1)$-module, which means that for all $n\in\Z$, the subspace of homogeneous degree-$n$ elements is finite-dimensional.
    
    We will show below in \cref{knn_monoidal} that $H_\ko^*$ is a symmetric monoidal functor via the \term{cross product} map $H_\ko^*(M)\otimes H_\ko^*(N)\to H_\ko^*(M\wedge_\ko N)$.
    %\TODO: symmetric monoidality. I think the K\"unneth map  the structure of a lax symmetric monoidal functor (but: reference?). The collapse of the K\"unneth spectral sequence implies all of those maps are isomorphisms -- does that suffice to show $H_\ko^*$ is strong symmetric monoidal? Are there any weird $\infty$-cateogry shenanigans here?
    \item Finally, $L_1$ is the $Q_1$-localization functor from \cref{local_stmod}, which Ricka~\cite[Proposition 2.6]{ricka_local_2016} showed is symmetric monoidal. This needs no finiteness assumptions, so we implicitly compose with the forgetful functor from finite type $\cA(1)$-modules to all $\cA(1)$-modules, which is symmetric monoidal.
\end{enumerate}
Thus $\Phi$ is a symmetric monoidal functor of $\infty$-categories, once we fulfill our promise about the cross product.

\begin{defn}
Let $\mu\colon H\Z/2\wedge_\ko H\Z/2\to H\Z/2$ be the multiplication map in the $E_\infty$-$\ko$-algebra structure on $H\Z/2$ we discussed above. For $\ko$-module maps $\phi_1\colon M\to \Sigma^{t_1} H\Z/2$ and $\phi_2\colon N\to\Sigma^{t_2} H\Z/2$, define their \term{cross product} $t_1\times t_2$ to be the composition
\begin{subequations}
\begin{equation}
    M\wedge_\ko N \xrightarrow{\phi_1\wedge\phi_2} \Sigma^{t_1} H\Z/2\wedge_\ko \Sigma^{t_2}H\Z/2\xrightarrow\simeq \Sigma^{t_1+t_2} H\Z/2\wedge_\ko H\Z/2\xrightarrow\mu \Sigma^{t_1+t_2}H\Z/2.
\end{equation}
Passing to homotopy groups of mapping spectra, we get a map
\begin{equation}
\label{cross}
\times\colon H_\ko^*(M)\otimes H_\ko^*(N)\to H_\ko^*(M\wedge_\ko N),
\end{equation}
\end{subequations}
which we also call the cross product.
\end{defn}

\begin{lem}
\label{knn_monoidal}
The functor $H_\ko^*\colon\Mod_\ko^{\mathit{ft}}\to\Mod_{\cA(1)}^{\mathit{ft}}$ is symmetric monoidal via the cross product $\times\colon H_\ko^*(M)\otimes H_\ko^*(N)\to H_\ko^*(M\wedge_\ko N)$. In particular, the cross product~\eqref{cross} is an isomorphism for all finite type $\ko$-modules $M$ and $N$.
\end{lem}

The proof of \cref{knn_monoidal} is identical to the usual proof of the standard fact that, by the Künneth formula, mod $2$ cohomology is a symmetric monoidal functor out of the $\infty$-category of finite type spaces. Absent the finite type hypothesis, $H^*(\bl;\Z/2)$ and $H_\ko^*$ are only lax symmetric monoidal.

%\begin{proof}
%\TODO: proof sketch, quite sketchy. AD (note to self): see meeting notes from 9/23. we can express this in terms of the universal coeffs sseq in EKMM, using that things are finite type.
%
%Since $H_\ko^*$ passes to homotopy classes of maps, it by definition factors through the homotopy category of $\Mod_\ko$. Passing to the homotopy category is symmetric monoidal (but, \TODO: find a reference). Now we just need to show that $H_\ko^*\colon\mathrm{Ho}(\Mod_\ko^{\mathit{ft}})\to\Mod_{\cA(1)}^{\mathit{ft}}$ is a symmetric monoidal functor -- much better, as its domain and codomain are $1$-categories.
%
%Anyways, the cross product map straighforwardly satisfies the axioms of a lax symmetric monoidal functor. For it to be symmetric monoidal, it needs to be  an equivalence. For this, use the K\"unneth spectral sequence, which instantly collapses since its $E^2$-page is Tor over $\Z/2$. % (\TODO: do we need to restrict to something like finite type $\cA(1)$-modules?)
%\end{proof}

\begin{lem}
\label{Phi_on_1001}
In $\StMod_{\cA(1)}[Q_1^{-1}]$, $\Phi(E_{1,0})\cong \Sigma T_1$ and $\Phi(E_{0,1})\cong \Sigma T_0$.
\end{lem}
\begin{proof}
After unwinding the definition, $\Phi(E_{1,0})$ is (the $Q_1$-localization of) $H_\ko^*$ of the Thom spectrum of the map $T\circ (x,0)\colon B\Z/2\to B\O/B\Spin\to B\GL_1(\ko)$, where $x\in H^1(B\Z/2;\Z/2)$ is the generator. If $\sigma\to B\Z/2$ denotes the tautological line bundle, then $(x,0) = (w_1(\sigma), w_2(\sigma))$, so $M\mathbf T(E_{1,0})\simeq \ko\wedge (B\Z/2)^{\sigma-1}$, as follows by combining~\cite[Example 1.6 and Lemma 1.8]{DY23}. For any spectrum $X$, there is a natural $\cA(1)$-module isomorphism $H_\ko^*(X)\cong H^*(X;\Z/2)$~\cite[Lemma 5.16]{BottSpiral}, where the $\cA(1)$-module structure on the latter is the restriction of the usual $\cA$-module structure. Finally, Campbell~\cite[Example 6.4]{Cam17} shows $H^*((B\Z/2)^{\sigma-1});\Z/2)\cong\Sigma T_1$. As this is already $Q_1$-local by~\cite[Chapter 2]{Yu95}, $L_1(\Sigma T_1)\cong \Sigma T_1$. As we saw in \cref{Q1_torsion}, $T_1$ is an invertible object of $\StMod_{\cA(1)}[Q_1^{-1}]$, and invertibility is closed under suspension (since $(\Sigma M)^{-1} = \Sigma^{-1}M^{-1}$).

For $E_{0,1}$, the story is very similar, but with $(x,x^2)$ in place of $(x,0)$ as the map to $B\O/B\Spin$, $-\sigma$ in place of $\sigma$, and $(B\Z/2)^{1-\sigma}$ in place of $(B\Z/2)^{\sigma-1}$. Campbell~\cite[Example 6.6]{Cam17} shows $H^*((B\Z/2)^{1-\sigma};\Z/2)\cong\Sigma T_0$, and  $\Sigma T_0$ is $Q_1$-local~\cite[Chapter 2]{Yu95}, so the $Q_1$-localization of $\Sigma T_0$ is itself, and it is an invertible object of $\StMod_{\cA(1)}[Q_1^{-1}]$, as displayed in \cref{Q1_torsion}.
\end{proof}

\begin{cor}
For all $\ell,k\in\N$, $\Phi(E_{\ell,k})$ is invertible.
\end{cor}
\begin{proof}
By \cref{Phi_on_1001}, $\Phi(E_{1,0})\cong \Sigma T_1$ and $\Phi(E_{0,1})\cong\Sigma T_0$, and by \cref{yu_thm}, both of these are invertible in $\StMod_{\cA(1)}[Q_1^{-1}]$. Since $\Phi$ is symmetric monoidal, $\Phi(E_{\ell,k})$ is a tensor product of copies of $\Phi(E_{1,0})$ and $\Phi(E_{0,1})$, and a tensor product of invertible objects is invertible.
\end{proof}
\begin{cor}
\label{our_phi_defn}
$\Phi$ restricts to a homomorphism of commutative monoids
\begin{equation}
    \Phi\colon\EMon \longrightarrow \Pic(\StMod_{\cA(1)}[Q_1^{-1}]). % \xrightarrow[\eqref{yu_thm}]{y} \Z\oplus\Z/4.
\end{equation}
\end{cor}
We will frequently compose $\Phi$ with the isomorphism $y\colon \Pic(\StMod_{\cA(1)}[Q_1^{-1}])\xrightarrow\cong\Z\oplus\Z/4$ from \cref{yu_thm}.

%%%%%%%%%%%%%%%%%%%%%
%%%%%%%%%%%%%%%%%%%%%
%\TODO: should we combine this and the next section into one section?

\section{Proofs of the main theorems}
\label{s:proofs}
% the first one: Pic(Q1) is the group completiom of the monoid of E_{ell,k} groups wrt fermionic tensor product
% Now we can prove o
Our first main result, \cref{Q1mainthm}, is that the map $y\circ \Phi\colon \EMon\to \Z\oplus\Z/4$ from \cref{our_phi_defn} is a group completion. That is, we will show that any map $\EMon$ to a group $G$ factors uniquely through $\Phi$.

%\Arun{TODO: make sure $E_{1,0}$ and $E_{0,1}$ map to the correct things}
\begin{proof}[Proof of \cref{Q1mainthm}]
Using the monoid isomorphism %$y\colon\Pic(\StMod_{\cA(1)}[Q_1^{-1}])\to \Z\oplus\Z/4$ and 
$\widetilde\varepsilon\colon \N^2/((4,0) - (0,4))\to\EMon$ from \cref{M_is_mod_4004}, it suffices to show that the following map is a group completion:
%we may replace $\Phi$ with the map
\begin{equation}
\label{just_numbers}
    \begin{aligned}
        y\circ \Phi\circ \widetilde\varepsilon\colon \N^2/((4,0) - (0,4)) &\longrightarrow \Z\oplus\Z/4\\
        (1,0) &\longmapsto (1,1)\\
        (0,1) &\longmapsto (1,0).
    \end{aligned}
\end{equation}
Let $G$ be a group and $f\colon\N^2/((4,0)-(0,4))\to G$ be a monoid homomorphism. We may without loss of generality assume $G$ is abelian and write elements of $G$ additively. Define
\begin{equation}
    \begin{aligned}
        \widetilde f\colon \Z\oplus\Z/4 &\longrightarrow G\\
        (a,b)&\longmapsto (a-b)f(0,1) + b f(1,0).
    \end{aligned}
\end{equation}
Since $f(0,4) = f(4,0)$ by the relation in the domain of $f$, $\widetilde f(a,b)$ only depends on $b\bmod 4$, so indeed $\widetilde f$ is well-defined as a homomorphism of abelian groups.

Since $\widetilde f(1,0) = f(0,1)$ and $\widetilde f(1,1) = f(1,0)$, $f = \widetilde f\circ y\circ \Phi\circ \widetilde\varepsilon$. To show this factorization is unique, suppose $\widetilde f'$ is another map with $\widetilde f'\circ y\circ \Phi\circ \widetilde\varepsilon = f$. Thus $\widetilde f'(1,0) = f(0,1) = \widetilde f(1,0)$ and $\widetilde f'(1,1) = f(1,0) = \widetilde f(1,1)$. That is, $\widetilde f$ and $\widetilde f'$ agree on a generating set $\{(1,0), (1,1)\}$ of $\Z\oplus\Z/4$, so $\widetilde f = \widetilde f'$.
\end{proof}

%%%%%%%%%%%%%%%%%%%%%
%%%%%%%%%%%%%%%%%%%%%

%\section{The diagonal subgroups and Smith periodicity}
\label{s:diag}

Finally, we compare the two $\Z/4$s that we introduced in \S\ref{ss:first_Z4} and \S\ref{ss:second_Z4}, which gave rise to our two periodicities. We do this in two ways. First, in \cref{calc_psi}, we compute the ``obvious'' map $\Z/4\cong \Pic(\StMod_{\cA(1)}[Q_1^{-1}])_{\mathrm{tors}}\to [B\Z/2, B\O/B\Spin]$, namely the map induced from $\delta\colon\EMon\to [B\Z/2, B\O/B\Spin]$ by the universal property of the group completion. This map turns out to be $2\colon\Z/4\to\Z/4$---not an isomorphism. After discussing what went wrong, we prove \cref{period_compare}, offering a different way to identify these two $\Z/4$s as generated by the fermionic group $E_{1,0}$.

%we prove \cref{period_compare}, identifying the two $\Z/4$s we introduced in .
\begin{prop}
\label{calc_psi}
The map $\Psi\colon\Z/4\to\Z/4$ defined in~\eqref{naive_comparison} is multiplication by $2$.
\end{prop}
\begin{proof}
Recall that $t\colon\Z/4\hookrightarrow \Z\oplus\Z/4$ is the inclusion of the torsion subgroup and $\widetilde\delta\colon \Z\oplus\Z/4\to [B\Z/2, B\O/B\Spin]$ is the unique map such that $\delta = \widetilde\delta\circ y \circ \Phi$, as guaranteed by \cref{Q1mainthm}.
We calculate $\widetilde\delta\circ t$ on the generator of $\Z/4$, which is the module $T_1$ by \cref{yu_thm}:
we have $\Phi(E_{1,0}) - \Phi(E_{0,1}) = T_1$, so 
\begin{equation}
    \widetilde\delta(0, 1) = \delta(E_{1,0}) - \delta(E_{0,1}).
\end{equation}
By \cref{rank1_delta}, $\delta(E_{1,0}) = (x, 0)$ and $\delta(E_{1,1}) = (x, x^2)$. As we saw in \cref{explicit_twists}, these correspond to $1$, resp.\ $3$ in $\Z/4$ under \cref{second_Z4}, so we conclude that $\widetilde\delta(0, 1) = 1 - 3 = 2\in\Z/4$.
\end{proof}
%So we have left to compute $\delta(E_{1,0})$ and $\delta(E_{0,1})$.
%
%
%
%
%
% We want to show that $\widetilde\delta\circ t\colon \Z/4\to [B\Z/2, B\O/B\Spin]$ sends a generator to twice a generator
%is an isomorphism. (By \cref{second_Z4}, we know the domain and codomain are isomorphic.) Here $t\colon\Z/4\hookrightarrow \Z\oplus\Z/4$ is the inclusion of the torsion subgroup and $\widetilde\delta\colon \Z\oplus\Z/4\to [B\O/B\Spin]$ is the unique map such that $\delta = \widetilde\delta\circ\Phi$, as guaranteed by \cref{Q1mainthm}. We will prove this by calculating $\widetilde\delta\circ t$ on the generator $1$ of $\Z/4$.
%
%The map $t$ is the inclusion of the torsion subgroup of $\Z\oplus\Z/4$, so $\widetilde\delta\circ t(1) = \widetilde\delta(0,1)$. Since $\widetilde\delta \circ\Phi = \delta$ and $(0,1) = \Phi(E_{1,0}) - \Phi(E_{0,1})$~\eqref{just_numbers}, then in $[B\Z/2, B\O/B\Spin]$,
%
%\TODO: this might not be a %generator?
%\TODO\dots

% \Cref{calc_psi} is not the result you'd expect knowing that these two $4$-periodicities are related; the canonical map ought to be an isomorphism, right?
One might have expected \cref{calc_psi} to be an isomorphism.
The issue stems from the different requirements for normalization between the Picard groupoid and the group of twists.
In the former, torsion elements might not begin in degree zero, but in the latter, twists are given by rank zero virtual bundles.
In particular,
% The discrepancy has to do with a degree shift---
by \cref{Phi_on_1001}, $\Phi(E_{1,0})$ is not a torsion element of $\Pic(\StMod_{\cA(1)}[Q_1^{-1}])$.
Instead, it is the \emph{suspension} of a torsion element. Generally, one does not think of $M$ and $\Sigma M$ as vastly different modules. But this means that $(0, 1)\in\Z\oplus\Z/4$ is not represented just by $E_{1,0}$: now, the suspension matters, and we have to represent $(0,1)$ by the formal difference $E_{1,0} - E_{0,1}$.

One would want to modify $\Psi$ to fix this, %in that the degree shifts
but the setup is very rigid: there is a unique $\Z/4$ subgroup of $\Pic(\StMod_{\cA(1)}[Q_1^{-1}])$, and the group completion is essentially unique. So, while $\Psi$ does not reflect how suspension shouldn't change the $\Z/4$ periodicity of $Q_1$-locally invertible $\cA(1)$-modules, there is not a good way to modify $\Psi$.
This is why we formulate \cref{period_compare} as a correspondence.

%\TODO: be clearer about the specific isomorphism $\Pic(\dots)\cong\Z\oplus\Z/4$. Like, are we writing $\widetilde{\delta}(\Phi(E_{\ell,k}))$ or just $\widetilde{\delta}(\ell,k)$, or what? Note 10/8 by AD: the codomain of the group completion map $\Phi$ is $\Z\oplus\Z/4$, so the domain of $\widetilde\delta$ is $\Z\oplus\Z/4$. Thus both $\widetilde\delta(\Phi(E_{1,0}))$ and $\widetilde\delta(\ell, k)$ parse. But maybe there's a less confusing way to write this

\begin{proof}[Proof of \cref{period_compare}]

We compute each side of the diagram.
On the left, we have $\widetilde\delta(y(\Phi(E_{1,0}))) = \delta(E_{1,0})$ by definition of $\widetilde\delta$. By \cref{rank1_delta}, $\delta(E_{1,0}) = (x, 0)$, and by \cref{explicit_twists}, this is a generator of $[B\Z/2, B\O/B\Spin]$.

On the right, by definition $\varsigma$ suspends or desuspends a module until it is stably equivalent to $T_0$, $T_1$, $T_2$, or $T_3$ (see \cref{desuspend}). By \cref{Phi_on_1001}, $\Phi(E_{1,0}) \cong \Sigma T_1$, so $\varsigma(\Phi(E_{1,0})) = T_1$, which by \cref{yu_thm} generates the torsion subgroup of $\Pic(\StMod_{\cA(1)}[Q_1^{-1}])$.
\end{proof}

%\TODO: Identify the $\Z/4$ Smith periodicity with the $\Z/4$ torsion summand of $\Pic(\cA(1)_{(Q_1)})$. Specifically, taking the diagonal subgroup defines a monoid homomorphism from $\EMon$ to $[B\Z/2, B\O/B\Spin]\cong\Z/4$. As the codomain is an abelian group, this factors through the group completion, which by \cref{Q1mainthm} gives us a map $\Pic(\cA(1)_{(Q_1)})\to [B\Z/2, B\O/B\Spin]$. Here we want to prove that restricted to the torsion subgroup of the Picard group, this map is an isomorphism $\Z/4\to\Z/4$, indeed giving a canonical isomorphism between these two periodicities.

%%%%%%%%%%%%%%%%%%%%%
%%%%%%%%%%%%%%%%%%%%%
   
\bibliographystyle{alpha}
\bibliography{CK_Zotero}

@article{freed_reflection_2021,
	title = {Reflection positivity and invertible topological phases},
	volume = {25},
	issn = {1364-0380, 1465-3060},
	url = {https://arxiv.org/abs/1604.06527},
	doi = {10.2140/gt.2021.25.1165},
	language = {en},
	number = {3},
	urldate = {2023-01-22},
	journal = {Geometry \& Topology},
	author = {Freed, Daniel S. and Hopkins, Michael J.},
	month = may,
	year = {2021},
	note = {\url{https://arxiv.org/abs/1604.06527}},
	pages = {1165--1330},
}

@misc{ablamowicz2016cliffordalgebrasrelatedfinite,
      title={On {C}lifford Algebras and Related Finite Groups and Group Algebras}, 
      author={Rafa\l{} Ab{\l}amowicz},
      year={2016},
      eprint={1610.03583},
      archivePrefix={arXiv},
      primaryClass={math.RA},
      note={\url{https://arxiv.org/abs/1610.03583}}, 
}

@article{albuquerque2002clifford,
  title={Clifford algebras obtained by twisting of group algebras},
  author={Albuquerque, Helena and Majid, Shahn},
  journal={Journal of Pure and Applied Algebra},
  volume={171},
  number={2-3},
  pages={133--148},
  year={2002},
  publisher={Elsevier},
    note = {\url{https://arxiv.org/abs/math/0011040}}
}

@inproceedings{kitaev_periodic_2009,
	title = {Periodic table for topological insulators and superconductors},
	url = {https://arxiv.org/abs/0901.2686},
	doi = {10.1063/1.3149495},
	language = {en},
	urldate = {2023-01-24},
	booktitle = {{AIP} {Conference} {Proceedings}},
	author = {Kitaev, Alexei},
	year = {2009},
	note = {\url{https://arxiv.org/abs/0901.2686}},
	pages = {22--30},
}

@article{atiyah_clifford_1963,
	AUTHOR = {Atiyah, M. F. and Bott, R. and Shapiro, A.},
     TITLE = {Clifford modules},
   JOURNAL = {Topology},
  FJOURNAL = {Topology. An International Journal of Mathematics},
    VOLUME = {3},
      YEAR = {1964},
    NUMBER = {suppl, suppl. 1},
     PAGES = {3--38},
      ISSN = {0040-9383},
   MRCLASS = {57.30 (55.30)},
  MRNUMBER = {167985},
MRREVIEWER = {J. F. Adams},
       DOI = {10.1016/0040-9383(64)90003-5},
       URL = {https://doi.org/10.1016/0040-9383(64)90003-5},
}

@article{freed_loop_2011,
	title = {Loop groups and twisted {K}-theory {I}},
	volume = {4},
	issn = {17538416},
	url = {https://arxiv.org/abs/0711.1906},
	doi = {10.1112/jtopol/jtr019},
	language = {en},
	number = {4},
	urldate = {2023-04-05},
	journal = {Journal of Topology},
	author = {Freed, Daniel S. and Hopkins, Michael J. and Teleman, Constantin},
	month = dec,
	year = {2011},
	note = {\url{https://arxiv.org/abs/0711.1906}},
	pages = {737--798},
}

@article{beardsley2023brauer,
 AUTHOR = {Beardsley, Jonathan and Luecke, Kiran and Morava, Jack},
     TITLE = {Brauer--{W}all groups and truncated {P}icard spectra of
              {K}-theory},
   JOURNAL = {Algebr. Geom. Topol.},
  FJOURNAL = {Algebraic \& Geometric Topology},
    VOLUME = {26},
      YEAR = {2026},
    NUMBER = {5},
     PAGES = {1749--1780},
      ISSN = {1472-2747,1472-2739},
   MRCLASS = {19L50 (55N15 55P42)},
  MRNUMBER = {5079112},
       DOI = {10.2140/agt.2026.26.1749},
       URL = {https://doi.org/10.2140/agt.2026.26.1749},
  note = {\url{https://arxiv.org/abs/2306.10112}},
}

@article{stehouwer_interacting_2022,
	title = {Interacting {SPT} phases are not {Morita} invariant},
	volume = {112},
	issn = {0377-9017, 1573-0530},
	url = {https://arxiv.org/abs/2110.07408},
	doi = {10.1007/s11005-022-01559-6},
	language = {en},
	number = {3},
	urldate = {2023-04-05},
	journal = {Letters in Mathematical Physics},
	author = {Stehouwer, Luuk},
	month = jun,
	year = {2022},
	note = {\url{https://arxiv.org/abs/2110.07408}},
	pages = {64},
}

@article{queiroz_dimensional_2016,
	title = {Dimensional hierarchy of fermionic interacting topological phases},
	volume = {117},
	issn = {0031-9007, 1079-7114},
	url = {https://arxiv.org/abs/1601.01596},
	doi = {10.1103/PhysRevLett.117.206405},
	language = {en},
	number = {20},
	urldate = {2024-04-17},
	journal = {Physical Review Letters},
	author = {Queiroz, Raquel and Khalaf, Eslam and Stern, Ady},
	month = nov,
	year = {2016},
	note = {\url{https://arxiv.org/abs/1601.01596}},
	pages = {206405},
}

@article{Thiang_2015,
   title={On the {K}-Theoretic Classification of Topological Phases of Matter},
   volume={17},
   ISSN={1424-0661},
   url={https://dx.doi.org/10.1007/s00023-015-0418-9},
   DOI={10.1007/s00023-015-0418-9},
   number={4},
   journal={Annales Henri Poincaré},
   publisher={Springer Science and Business Media LLC},
   author={Thiang, Guo Chuan},
   year={2015},
   month=may, pages={757–794},
    note = {\url{https://arxiv.org/abs/1406.7366}}
}

@article{debray_long_2024,
	 AUTHOR = {Debray, Arun and Devalapurkar, Sanath K. and Krulewski,
              Cameron and Liu, Yu Leon and Pacheco-Tallaj, Natalia and
              Thorngren, Ryan},
     TITLE = {A long exact sequence in symmetry breaking: order parameter
              constraints, defect anomaly-matching, and higher {B}erry
              phases},
   JOURNAL = {J. High Energy Phys.},
  FJOURNAL = {Journal of High Energy Physics},
      YEAR = {2025},
    NUMBER = {7},
     PAGES = {Paper No. 7, 47},
      ISSN = {1126-6708,1029-8479},
   MRCLASS = {81T50 (81Q70 81T45)},
  MRNUMBER = {4927901},
       DOI = {10.1007/jhep07(2025)007},
	note = {\url{https://arxiv.org/abs/2309.16749}},
}

@article{debray_smith_2024,
	AUTHOR = {Debray, Arun and Devalapurkar, Sanath K. and Krulewski,
              Cameron and Liu, Yu Leon and Pacheco-Tallaj, Natalia and
              Thorngren, Ryan},
     TITLE = {The {S}mith Fiber Sequence and Invertible Field
              Theories},
   JOURNAL = {Comm. Math. Phys.},
  FJOURNAL = {Communications in Mathematical Physics},
    VOLUME = {407},
      YEAR = {2026},
    NUMBER = {2},
     PAGES = {Paper No. 25},
      ISSN = {0010-3616,1432-0916},
   MRCLASS = {57 (81T45)},
  MRNUMBER = {5015202},
       DOI = {10.1007/s00220-025-05505-0},
	note = {\url{https://arxiv.org/abs/2405.04649}},
}

@incollection{Bru14,
    AUTHOR = {Bruner, Robert R.},
     TITLE = {Idempotents, localizations and {P}icard groups of {$
              A(1)$}-modules},
 BOOKTITLE = {An alpine expedition through algebraic topology},
    SERIES = {Contemp. Math.},
    VOLUME = {617},
     PAGES = {81--108},
 PUBLISHER = {Amer. Math. Soc., Providence, RI},
      YEAR = {2014},
   MRCLASS = {55S10 (19L41)},
  MRNUMBER = {3243394},
MRREVIEWER = {Constanze Roitzheim},
       DOI = {10.1090/conm/617/12346},
	note = {\url{https://arxiv.org/abs/1211.0213}}
}

@book{BG10,
    AUTHOR = {Bruner, Robert R. and Greenlees, J. P. C.},
     TITLE = {Connective real {$K$}-theory of finite groups},
    SERIES = {Mathematical Surveys and Monographs},
    VOLUME = {169},
 PUBLISHER = {American Mathematical Society, Providence, RI},
      YEAR = {2010},
     PAGES = {vi+318},
      ISBN = {978-0-8218-5189-0},
   MRCLASS = {19L41 (13D45 19-02 19L47 19L64 55N15 55N91)},
  MRNUMBER = {2723113},
       DOI = {10.1090/surv/169},
       URL = {https://doi.org/10.1090/surv/169},
}

@article{BL01,
    AUTHOR = {Baker, Andrew and Lazarev, Andrej},
     TITLE = {On the {A}dams spectral sequence for {$R$}-modules},
   JOURNAL = {Algebr. Geom. Topol.},
  FJOURNAL = {Algebraic \& Geometric Topology},
    VOLUME = {1},
      YEAR = {2001},
     PAGES = {173--199},
      ISSN = {1472-2747},
   MRCLASS = {55P42 (55N20 55P43 55T15)},
  MRNUMBER = {1823498},
MRREVIEWER = {Mark Hovey},
       DOI = {10.2140/agt.2001.1.173},
        note = "\url{https://arxiv.org/abs/math/0105079}"
}

@article{Bak20,
        author = {Andrew Baker},
        title = {Homotopy theory of modules over a commutative {$S$}-algebra: some tools
  and examples},
        year = {2020},
        note = {\url{https://arxiv.org/abs/2003.12003}}
}

@article{Sto88,
    AUTHOR = {Stolz, Stephan},
     TITLE = {Exotic structures on {$4$}-manifolds detected by spectral
              invariants},
   JOURNAL = {Invent. Math.},
  FJOURNAL = {Inventiones Mathematicae},
    VOLUME = {94},
      YEAR = {1988},
    NUMBER = {1},
     PAGES = {147--162},
      ISSN = {0020-9910},
   MRCLASS = {57R55 (57R85 58G25)},
  MRNUMBER = {958593},
MRREVIEWER = {J\"{u}rgen Eichhorn},
       DOI = {10.1007/BF01394348},
       URL = {https://doi.org/10.1007/BF01394348},
}

@article{Cam17,
        author = {Jonathan A. Campbell},
        title = {Homotopy Theoretic Classification of Symmetry Protected Phases},
        year = {2017},
        note = {\url{https://arxiv.org/abs/1708.04264}}
}

@article{AM71,
    AUTHOR = {Adams, J. F. and Margolis, H. R.},
     TITLE = {Modules over the {S}teenrod algebra},
   JOURNAL = {Topology},
  FJOURNAL = {Topology. An International Journal of Mathematics},
    VOLUME = {10},
      YEAR = {1971},
     PAGES = {271--282},
      ISSN = {0040-9383},
   MRCLASS = {18G40 (55G10)},
  MRNUMBER = {294450},
MRREVIEWER = {A. K. Bousfield},
       DOI = {10.1016/0040-9383(71)90020-6},
       URL = {https://doi.org/10.1016/0040-9383(71)90020-6},
}

@phdthesis{Yu95,
    AUTHOR = {Yu, Cherng-Yih},
     TITLE = {The connective real {K}-theory of elementary abelian 2-groups},
      YEAR = {1995},
    school = {The University of Notre Dame}
}

@book{Mar83,
    AUTHOR = {Margolis, H. R.},
     TITLE = {Spectra and the {S}teenrod algebra},
    SERIES = {North-Holland Mathematical Library},
    VOLUME = {29},
      NOTE = {Modules over the Steenrod algebra and the stable homotopy
              category},
 PUBLISHER = {North-Holland Publishing Co., Amsterdam},
      YEAR = {1983},
     PAGES = {xix+489},
      ISBN = {0-444-86516-0},
   MRCLASS = {55-02 (55Sxx)},
  MRNUMBER = {738973},
MRREVIEWER = {Andrzej Kozlowski},
}

@article{DYY23,
	AUTHOR = {Debray, Arun and Ye, Weicheng and Yu, Matthew},
     TITLE = {Bosonization and anomaly indicators of {$(2+1)$}-{D} fermionic
              topological orders},
   JOURNAL = {Comm. Math. Phys.},
  FJOURNAL = {Communications in Mathematical Physics},
    VOLUME = {406},
      YEAR = {2025},
    NUMBER = {8},
     PAGES = {Paper No. 178, 88},
      ISSN = {0010-3616,1432-0916},
   MRCLASS = {81T45 (18M20 57K16 82D03)},
  MRNUMBER = {4927807},
       DOI = {10.1007/s00220-025-05344-z},
       URL = {https://doi.org/10.1007/s00220-025-05344-z},
	note = {\url{https://arxiv.org/abs/2312.13341}}
}

@article{DYY,
		author = {Arun Debray and Weicheng Ye and Matthew Yu},
	title = {Global Structure in the Presence of a Topological Defect},
	JOURNAL = {Comm. Math. Phys.},
  FJOURNAL = {Communications in Mathematical Physics},
    VOLUME = {407},
      YEAR = {2026},
    NUMBER = {8},
     PAGES = {Paper No. 159},
      ISSN = {0010-3616,1432-0916},
   MRCLASS = {81T45 (57K16 82B10)},
  MRNUMBER = {5099898},
       DOI = {10.1007/s00220-026-05629-x},
       URL = {https://doi.org/10.1007/s00220-026-05629-x},
	note = {\url{https://arxiv.org/abs/2501.18399}}
}

@article{ABG18,
    AUTHOR = {Ando, Matthew and Blumberg, Andrew J. and Gepner, David},
     TITLE = {Parametrized spectra, multiplicative {T}hom spectra and the
              twisted {U}mkehr map},
   JOURNAL = {Geom. Topol.},
  FJOURNAL = {Geometry \& Topology},
    VOLUME = {22},
      YEAR = {2018},
    NUMBER = {7},
     PAGES = {3761--3825},
      ISSN = {1465-3060},
   MRCLASS = {55P43 (55R70)},
  MRNUMBER = {3890766},
MRREVIEWER = {Birgit Richter},
       DOI = {10.2140/gt.2018.22.3761},
       note = {\url{https://arxiv.org/abs/1112.2203}}
}

@article{Wan08,
    AUTHOR = {Wang, Bai-Ling},
     TITLE = {Geometric cycles, index theory and twisted {$K$}-homology},
   JOURNAL = {J. Noncommut. Geom.},
  FJOURNAL = {Journal of Noncommutative Geometry},
    VOLUME = {2},
      YEAR = {2008},
    NUMBER = {4},
     PAGES = {497--552},
      ISSN = {1661-6952},
   MRCLASS = {19K56 (19K33 55N22 58J22)},
  MRNUMBER = {2438341},
MRREVIEWER = {Thomas Schick},
       DOI = {10.4171/JNCG/27},
        note = {\url{https://arxiv.org/abs/0710.1625}}
}

@article{DY23,
        author = {Arun Debray and Matthew Yu},
        title = {Adams spectral sequences for non-vector-bundle {T}hom spectra},
        year = {2023},
        note = {\url{https://arxiv.org/abs/2305.01678}}
}

@book{MQRT77,
    AUTHOR = {May, J. Peter},
     TITLE = {{$E\sb{\infty }$} ring spaces and {$E\sb{\infty }$} ring
              spectra},
    SERIES = {Lecture Notes in Mathematics, Vol. 577},
      NOTE = {With contributions by Frank Quinn, Nigel Ray, and Jørgen
              Tornehave. \url{https://www.math.uchicago.edu/~may/BOOKS/e_infty.pdf}},
 PUBLISHER = {Springer-Verlag, Berlin-New York},
      YEAR = {1977},
     PAGES = {268},
   MRCLASS = {55D35},
  MRNUMBER = {0494077},
MRREVIEWER = {J. Stasheff},
}

@book{Pet68,
  author    = {F. P. Peterson},
  title     = {Lectures on Cobordism Theory},
  publisher = {Kinokuniya Book Store Co., Ltd.},
  year      = {1968},
  series    = {Lectures in Mathematics},
}

@incollection{Joa04,
    AUTHOR = {Joachim, Michael},
     TITLE = {Higher coherences for equivariant {$K$}-theory},
 BOOKTITLE = {Structured ring spectra},
    SERIES = {London Math. Soc. Lecture Note Ser.},
    VOLUME = {315},
     PAGES = {87--114},
 PUBLISHER = {Cambridge Univ. Press, Cambridge},
      YEAR = {2004},
   MRCLASS = {19L47 (55P43)},
  MRNUMBER = {2122155},
MRREVIEWER = {Kazuhisa Shimakawa},
       DOI = {10.1017/CBO9780511529955.006},
       URL = {https://doi.org/10.1017/CBO9780511529955.006},
}

@article{ABGHR14a,
    AUTHOR = {Ando, Matthew and Blumberg, Andrew J. and Gepner, David and
              Hopkins, Michael J. and Rezk, Charles},
     TITLE = {An {$\infty$}-categorical approach to {$R$}-line bundles,
              {$R$}-module {T}hom spectra, and twisted {$R$}-homology},
   JOURNAL = {J. Topol.},
  FJOURNAL = {Journal of Topology},
    VOLUME = {7},
      YEAR = {2014},
    NUMBER = {3},
     PAGES = {869--893},
      ISSN = {1753-8416},
   MRCLASS = {55P43 (55N20 55U40)},
  MRNUMBER = {3252967},
MRREVIEWER = {Tyler D. Lawson},
       DOI = {10.1112/jtopol/jtt035},
        note = "\url{https://arxiv.org/abs/1403.4325}"
}

@article{ABGHR14b,
    AUTHOR = {Ando, Matthew and Blumberg, Andrew J. and Gepner, David and
              Hopkins, Michael J. and Rezk, Charles},
     TITLE = {Units of ring spectra, orientations and {T}hom spectra via
              rigid infinite loop space theory},
   JOURNAL = {J. Topol.},
  FJOURNAL = {Journal of Topology},
    VOLUME = {7},
      YEAR = {2014},
    NUMBER = {4},
     PAGES = {1077--1117},
      ISSN = {1753-8416},
   MRCLASS = {55P43 (55N34 55R65)},
  MRNUMBER = {3286898},
MRREVIEWER = {Steven R. Costenoble},
       DOI = {10.1112/jtopol/jtu009},
        note = "\url{https://arxiv.org/abs/1403.4320}"
}

@misc{Sto98,
        author = {Stephan Stolz},
        title = {Concordance classes of positive scalar curvature metrics},
        year = {1998},
        note = {\url{https://www3.nd.edu/~stolz/preprint.html}}
}

@article{HJ20,
    AUTHOR = {Hebestreit, Fabian and Joachim, Michael},
     TITLE = {Twisted spin cobordism and positive scalar curvature},
   JOURNAL = {J. Topol.},
  FJOURNAL = {Journal of Topology},
    VOLUME = {13},
      YEAR = {2020},
    NUMBER = {1},
     PAGES = {1--58},
      ISSN = {1753-8416},
   MRCLASS = {55N22 (55R70)},
  MRNUMBER = {3999671},
MRREVIEWER = {Cenap \"{O}zel},
       DOI = {10.1112/topo.12122},
        note = "\url{https://arxiv.org/abs/1311.3164}"
}

@article{TY19,
        title={Why are fractional charges of orientifolds compatible with {D}irac quantization?},
        author={Yuji Tachikawa and Kazuya Yonekura},
        journal={SciPost Phys.},
        volume={7},
        issue={5},
        pages={58},
        year={2019},
        publisher={SciPost},
        doi={10.21468/SciPostPhys.7.5.058},
        note = "\url{https://arxiv.org/abs/1805.02772}"
}

@article{IT23,
    AUTHOR = {Ichikawa, Takumi and Tachikawa, Yuji},
     TITLE = {The super {F}robenius-{S}chur indicator and finite group gauge
              theories on {$\rm pin^-$} surfaces},
   JOURNAL = {Comm. Math. Phys.},
  FJOURNAL = {Communications in Mathematical Physics},
    VOLUME = {400},
      YEAR = {2023},
    NUMBER = {1},
     PAGES = {417--428},
      ISSN = {0010-3616},
   MRCLASS = {20C15 (16W55 20C35 81T13)},
  MRNUMBER = {4581479},
MRREVIEWER = {Alexander Isaakovich Shtern},
       DOI = {10.1007/s00220-022-04601-9},
       note = {\url{https://arxiv.org/abs/2002.10642}}
}

@article{Ben88,
    AUTHOR = {Benson, Dave},
     TITLE = {Spin modules for symmetric groups},
   JOURNAL = {J. London Math. Soc. (2)},
  FJOURNAL = {Journal of the London Mathematical Society. Second Series},
    VOLUME = {38},
      YEAR = {1988},
    NUMBER = {2},
     PAGES = {250--262},
      ISSN = {0024-6107},
   MRCLASS = {20C30},
  MRNUMBER = {966297},
MRREVIEWER = {A. Kerber},
       DOI = {10.1112/jlms/s2-38.2.250},
       URL = {https://doi.org/10.1112/jlms/s2-38.2.250},
}

@incollection{Boa99,
    AUTHOR = {Boardman, J. Michael},
     TITLE = {Conditionally convergent spectral sequences},
 BOOKTITLE = {Homotopy invariant algebraic structures ({B}altimore, {MD},
              1998)},
    SERIES = {Contemp. Math.},
    VOLUME = {239},
     PAGES = {49--84},
 PUBLISHER = {Amer. Math. Soc., Providence, RI},
      YEAR = {1999},
   MRCLASS = {55T05 (18A30 18G40)},
  MRNUMBER = {1718076},
MRREVIEWER = {John McCleary},
       DOI = {10.1090/conm/239/03597},
}

@article{COSY,
        doi = {10.1103/physrevd.102.045019}, 
        url = {https://doi.org/10.1103%2Fphysrevd.102.045019},
        year = 2020,
        month = {aug},
        publisher = {American Physical Society ({APS})},
        volume = {102},
        number = {4},
        author = {Clay Córdova and Kantaro Ohmori and Shu-Heng Shao and Fei Yan},
        title = {Decorated {$\Z_2$} symmetry defects and their time-reversal anomalies},
        journal = {Physical Review D},
        note = {\url{https://arxiv.org/abs/1910.14046}}
}

@article{HA,
note = {\url{https://www.math.ias.edu/~lurie/papers/HA.pdf}},
author = {Jacob Lurie},
title = {Higher Algebra},
year = {2017}
}

@article{karoubi_algebres_1968,
	AUTHOR = {Karoubi, Max},
     TITLE = {Algèbres de {C}lifford et {$K$}-théorie},
   JOURNAL = {Ann. Sci. École Norm. Sup. (4)},
  FJOURNAL = {Annales Scientifiques de l'\École Normale Supérieure. Quatrième
              Série},
    VOLUME = {1},
      YEAR = {1968},
     PAGES = {161--270},
      ISSN = {0012-9593},
   MRCLASS = {18.20 (55.00)},
  MRNUMBER = {238927},
MRREVIEWER = {J. F. Adams},
       URL = {https://www.numdam.org/item?id=ASENS_1968_4_1_2_161_0},
}

@article{AP76,
    AUTHOR = {Adams, J. F. and Priddy, S. B.},
     TITLE = {Uniqueness of {$B\SO$}},
   JOURNAL = {Math. Proc. Cambridge Philos. Soc.},
  FJOURNAL = {Mathematical Proceedings of the Cambridge Philosophical
              Society},
    VOLUME = {80},
      YEAR = {1976},
    NUMBER = {3},
     PAGES = {475--509},
      ISSN = {0305-0041},
   MRCLASS = {55D35},
  MRNUMBER = {431152},
MRREVIEWER = {V. P. Snaith},
       DOI = {10.1017/S0305004100053111},
       URL = {https://doi-org.ezproxy.lib.utexas.edu/10.1017/S0305004100053111},
}

@article{kreck_surgery_1999,
    AUTHOR = {Kreck, Matthias},
     TITLE = {Surgery and duality},
   JOURNAL = {Ann. of Math. (2)},
  FJOURNAL = {Annals of Mathematics. Second Series},
    VOLUME = {149},
      YEAR = {1999},
    NUMBER = {3},
     PAGES = {707--754},
      ISSN = {0003-486X},
   MRCLASS = {57R67 (57R65)},
  MRNUMBER = {1709301},
MRREVIEWER = {Laurence R. Taylor},
       DOI = {10.2307/121071},
       note = {\url{https://arxiv.org/abs/math/9905211}}
}

@article{madsen_infinite_1977,
    AUTHOR = {Madsen, I. and Snaith, V. and Tornehave, J.},
     TITLE = {Infinite loop maps in geometric topology},
   JOURNAL = {Math. Proc. Cambridge Philos. Soc.},
  FJOURNAL = {Mathematical Proceedings of the Cambridge Philosophical
              Society},
    VOLUME = {81},
      YEAR = {1977},
    NUMBER = {3},
     PAGES = {399--430},
      ISSN = {0305-0041},
   MRCLASS = {55D35 (55F45)},
  MRNUMBER = {494076},
MRREVIEWER = {J. P. May},
       DOI = {10.1017/S0305004100053482},
       URL = {https://doi.org/10.1017/S0305004100053482},
}

@article{antieau_actions_2014,
    AUTHOR = {Antieau, Benjamin and Gepner, David and Gómez, José Manuel},
     TITLE = {Actions of {$K(\pi,n)$} spaces on {$K$}-theory and uniqueness
              of twisted {$K$}-theory},
   JOURNAL = {Trans. Amer. Math. Soc.},
  FJOURNAL = {Transactions of the American Mathematical Society},
    VOLUME = {366},
      YEAR = {2014},
    NUMBER = {7},
     PAGES = {3631--3648},
      ISSN = {0002-9947},
   MRCLASS = {19L50 (55N15)},
  MRNUMBER = {3192610},
MRREVIEWER = {Bin Zhang},
       DOI = {10.1090/S0002-9947-2014-05937-0},
note = "\url{https://arxiv.org/abs/1106.5099}"
}

@article{lin_topological_2024,
	author = {Ying-Hsuan Lin and Mayuko Yamashita},
	title = {Topological Elliptic Genera {I} -- The mathematical foundation},
	year = {2024},
	note = {\url{https://arxiv.org/abs/2412.02298}}
}

@article{lastovecki_cohomology_2005,
    AUTHOR = {Lastovecki, Richard},
     TITLE = {Cohomology of {$B{\rm O}(n_1)\times\dots\times B{\rm O}(n_m)$}
              with local integer coefficients},
   JOURNAL = {Comment. Math. Univ. Carolin.},
  FJOURNAL = {Commentationes Mathematicae Universitatis Carolinae},
    VOLUME = {46},
      YEAR = {2005},
    NUMBER = {1},
     PAGES = {21--32},
      ISSN = {0010-2628,1213-7243},
   MRCLASS = {55R40 (55N25 55R45)},
  MRNUMBER = {2175856},
MRREVIEWER = {Andrew\ J.\ Baker},
}

@article{donovan_graded_1970,
    title = {Graded {B}rauer groups and {K}-theory with local coefficients},
    volume = {38},
    copyright = {http://www.springer.com/tdm},
    issn = {0073-8301, 1618-1913},
    url = {http://link.springer.com/10.1007/BF02684650},
    doi = {10.1007/BF02684650},
    language = {en},
    number = {1},
    urldate = {2024-07-02},
    journal = {Publications mathématiques de l'IHÉS},
    author = {Donovan, P. and Karoubi, M.},
    month = jan,
    year = {1970},
    pages = {5--25},
}

@article{kapustin_fermionic_2015,
    AUTHOR = {Kapustin, Anton and Thorngren, Ryan and Turzillo, Alex and
              Wang, Zitao},
     TITLE = {Fermionic symmetry protected topological phases and
              cobordisms},
   JOURNAL = {J. High Energy Phys.},
  FJOURNAL = {Journal of High Energy Physics},
      YEAR = {2015},
    NUMBER = {12},
     PAGES = {052, front matter+20pp},
      ISSN = {1126-6708,1029-8479},
   MRCLASS = {81T60},
  MRNUMBER = {3464750},
       DOI = {10.1007/jhep12(2015)052},
       note = {\url{https://arxiv.org/abs/1406.7329}}
}

@article{BottSpiral,
	author = {Arun Debray and Cameron Krulewski and Luuk Stehouwer},
	title = {Unraveling the {B}ott spiral},
	year = {2026},
	note = {\url{https://arxiv.org/abs/2605.00316}}
}

@article{wall_characterization_1962,
    AUTHOR = {Wall, C. T. C.},
     TITLE = {A characterization of simple modules over the {S}teenrod
              algebra {${\rm mod}\ 2$}},
   JOURNAL = {Topology},
  FJOURNAL = {Topology. An International Journal of Mathematics},
    VOLUME = {1},
      YEAR = {1962},
     PAGES = {249--254},
      ISSN = {0040-9383},
   MRCLASS = {55.34},
  MRNUMBER = {148061},
MRREVIEWER = {R.\ Brown},
       DOI = {10.1016/0040-9383(62)90108-8},
       URL = {https://doi.org/10.1016/0040-9383(62)90108-8},
}

@article{kuroda_computations_2025,
	author = {Naoki Kuroda},
	title = {Computations of {Spin-Sp(4)}, {Spin-SU(8)}, and {Spin-Spin(16)} bordism groups in dimensions up to 7},
	year = {2025},
	note = {\url{https://arxiv.org/abs/2504.15014}}
}

@article{gomez_nonexistence_2010,
    AUTHOR = {G{\'{o}}mez, José Manuel},
     TITLE = {On the nonexistence of higher twistings over a point for
              {B}orel cohomology of {$K$}-theory},
   JOURNAL = {J. Geom. Phys.},
  FJOURNAL = {Journal of Geometry and Physics},
    VOLUME = {60},
      YEAR = {2010},
    NUMBER = {4},
     PAGES = {678--685},
      ISSN = {0393-0440,1879-1662},
   MRCLASS = {19L50 (19L41 19L47 55N15 57T10)},
  MRNUMBER = {2602381},
MRREVIEWER = {Francis\ Clarke},
       DOI = {10.1016/j.geomphys.2009.12.009},
       URL = {https://doi.org/10.1016/j.geomphys.2009.12.009},
       note = {\url{https://arxiv.org/abs/0905.3151}}
}

@article{dadarlat_unit_2015,
    AUTHOR = {Dadarlat, Marius and Pennig, Ulrich},
     TITLE = {Unit spectra of {$K$}-theory from strongly self-absorbing
              {$C^*$}-algebras},
   JOURNAL = {Algebr. Geom. Topol.},
  FJOURNAL = {Algebraic \& Geometric Topology},
    VOLUME = {15},
      YEAR = {2015},
    NUMBER = {1},
     PAGES = {137--168},
      ISSN = {1472-2747,1472-2739},
   MRCLASS = {46L80 (19K99 55N15 55P42)},
  MRNUMBER = {3325734},
MRREVIEWER = {Valentin\ Deaconu},
       DOI = {10.2140/agt.2015.15.137},
       URL = {https://doi.org/10.2140/agt.2015.15.137},
       note = {\url{https://arxiv.org/abs/1306.2583}}
}

@article{pennig_noncommutative_2016,
    AUTHOR = {Pennig, Ulrich},
     TITLE = {A non-commutative model for higher twisted {$K$}-theory},
   JOURNAL = {J. Topol.},
  FJOURNAL = {Journal of Topology},
    VOLUME = {9},
      YEAR = {2016},
    NUMBER = {1},
     PAGES = {27--50},
      ISSN = {1753-8416,1753-8424},
   MRCLASS = {46L80 (19L50)},
  MRNUMBER = {3465839},
MRREVIEWER = {Magnus\ Goffeng},
       DOI = {10.1112/jtopol/jtv033},
       URL = {https://doi.org/10.1112/jtopol/jtv033},
       note = {\url{https://arxiv.org/abs/1502.02807}}
}

@mastersthesis{brook_higher_2020,
  author       = {David Brook},
  title = {Computations in higher twisted {$K$}-theory},
  school       = {The University of Adelaide},
  year   = {2020},
    note = {\url{https://arxiv.org/abs/2007.08964}}
}

@article{beaudry_guide_2024,
	 AUTHOR = {Beaudry, Agnès and Lewis, Chloe and May, Clover and Pauli,
              Sabrina and Tatum, Elizabeth},
     TITLE = {A guide to equivariant parametrized cohomology},
   JOURNAL = {Topology Appl.},
  FJOURNAL = {Topology and its Applications},
    VOLUME = {376},
      YEAR = {2025},
     PAGES = {Paper No. 109449, 71},
      ISSN = {0166-8641,1879-3207},
   MRCLASS = {55N91},
  MRNUMBER = {4965750},
       DOI = {10.1016/j.topol.2025.109449},
       URL = {https://doi.org/10.1016/j.topol.2025.109449},
	note = {\url{https://arxiv.org/abs/2410.13971}}
}

@article{beaudry_parametrized_2026,
	author = {Agnès Beaudry and Chloe Lewis and Clover May and Sabrina Pauli and Elizabeth Tatum},
	title = {Parametrized {$\underline{\mathbb{F}}_2$}-Cohomology of {$B_{C_2}O(1)$}},
	year = {2026},
	note = {\url{https://arxiv.org/abs/2607.27398}}
}

@article{cadek_cohomology_1999,
    AUTHOR = {{\v C}adek, Martin},
     TITLE = {The cohomology of {${\rm BO}(n)$} with twisted integer
              coefficients},
   JOURNAL = {J. Math. Kyoto Univ.},
  FJOURNAL = {Journal of Mathematics of Kyoto University},
    VOLUME = {39},
      YEAR = {1999},
    NUMBER = {2},
     PAGES = {277--286},
      ISSN = {0023-608X},
   MRCLASS = {55R40 (57R20)},
  MRNUMBER = {1709293},
MRREVIEWER = {Hans-Werner\ Henn},
       DOI = {10.1215/kjm/1250517912},
       URL = {https://doi.org/10.1215/kjm/1250517912},
}

@article{wendt_chow_2024,
    AUTHOR = {Wendt, Matthias},
     TITLE = {Chow-{W}itt rings of {G}rassmannians},
   JOURNAL = {Algebr. Geom. Topol.},
  FJOURNAL = {Algebraic \& Geometric Topology},
    VOLUME = {24},
      YEAR = {2024},
    NUMBER = {1},
     PAGES = {1--48},
      ISSN = {1472-2747,1472-2739},
   MRCLASS = {14C15 (14C17 14F43 14M15)},
  MRNUMBER = {4721362},
MRREVIEWER = {Aigli\ Papantonopoulou},
       DOI = {10.2140/agt.2024.24.1},
       URL = {https://doi.org/10.2140/agt.2024.24.1},
       note ={\url{https://arxiv.org/abs/1805.06142}}
}

@article{brazelton_chow_2025,
	author = {Thomas Brazelton and Matthias Wendt},
	title = {The {C}how--{W}itt rings of the classifying space of quadratically oriented bundles},
	year = {2025},
	note = {\url{https://arxiv.org/abs/2501.06302}}
}

@article{douglas_twisted_2006,
    AUTHOR = {Douglas, Christopher L.},
     TITLE = {On the twisted {$K$}-homology of simple {L}ie groups},
   JOURNAL = {Topology},
  FJOURNAL = {Topology. An International Journal of Mathematics},
    VOLUME = {45},
      YEAR = {2006},
    NUMBER = {6},
     PAGES = {955--988},
      ISSN = {0040-9383},
   MRCLASS = {19L41 (19D55 55N15 55N22)},
  MRNUMBER = {2263220},
MRREVIEWER = {Gregory\ D.\ Landweber},
       DOI = {10.1016/j.top.2006.06.007},
       URL = {https://doi.org/10.1016/j.top.2006.06.007},
       note = {\url{https://arxiv.org/abs/math/0402082}}
}

@article{costenoble_algebraic_2024,
    AUTHOR = {Costenoble, Steven R. and Hudson, Thomas and Tilson, Sean},
     TITLE = {An algebraic {$C_2$}-equivariant {B}\'{e}zout theorem},
   JOURNAL = {Algebr. Geom. Topol.},
  FJOURNAL = {Algebraic \& Geometric Topology},
    VOLUME = {24},
      YEAR = {2024},
    NUMBER = {4},
     PAGES = {2331--2350},
      ISSN = {1472-2747,1472-2739},
   MRCLASS = {55N91 (14N10 14N15 55R40 55R91)},
  MRNUMBER = {4776283},
MRREVIEWER = {Juan\ Antonio\ P\'{e}rez},
       DOI = {10.2140/agt.2024.24.2331},
       URL = {https://doi.org/10.2140/agt.2024.24.2331},
       note = {\url{https://arxiv.org/abs/2211.05382}}
}

@article{costenoble_graded_2023,
    AUTHOR = {Costenoble, Steven R.},
     TITLE = {The {$RO(\Pi B)$}-graded {$C_{2}$}-equivariant ordinary
              cohomology of {$B_{C_2}U(1)$}},
   JOURNAL = {Topology Appl.},
  FJOURNAL = {Topology and its Applications},
    VOLUME = {338},
      YEAR = {2023},
     PAGES = {Paper No. 108660, 53},
      ISSN = {0166-8641,1879-3207},
   MRCLASS = {55N91 (14N15 55R40 55R91)},
  MRNUMBER = {4629789},
MRREVIEWER = {Juan\ Antonio\ P\'{e}rez},
       DOI = {10.1016/j.topol.2023.108660},
       URL = {https://doi.org/10.1016/j.topol.2023.108660},
       note = {\url{https://arxiv.org/abs/1312.0926}}
}

@article{costenoble_equivariant_2025,
	author = {Steven R. Costenoble and Thomas Hudson},
	title = {The {$C_2$}-equivariant ordinary cohomology of complex quadrics {III}: Exceptional cases},
	year = {2025},
	note = {\url{https://arxiv.org/abs/2511.14264}}
}

@article{costenoble_geometric_2023,
	author = {Steven R. Costenoble and Thomas Hudson},
	title = {A geometric {$C_2$}-equivariant {B}ézout Theorem},
	year = {2023},
	note = {\url{https://arxiv.org/abs/2312.00559}}
}

@article{costenoble_equivariant_2022,
    AUTHOR = {Costenoble, Steven R. and Hudson, Thomas and Tilson, Sean},
     TITLE = {The {$C_2$}-equivariant cohomology of complex projective
              spaces},
   JOURNAL = {Adv. Math.},
  FJOURNAL = {Advances in Mathematics},
    VOLUME = {398},
      YEAR = {2022},
     PAGES = {Paper No. 108245, 69},
      ISSN = {0001-8708,1090-2082},
   MRCLASS = {55N91 (14N10 14N15 55R40 55R91)},
  MRNUMBER = {4388953},
MRREVIEWER = {No\'{e}\ B\'{a}rcenas Torres},
       DOI = {10.1016/j.aim.2022.108245},
       URL = {https://doi.org/10.1016/j.aim.2022.108245},
       note = {\url{https://arxiv.org/abs/1811.07355}}
}

@article{costenoble_ordinary_2017,
	author = {Steven R. Costenoble},
	title = {The {$\mathbb{Z}/p$} ordinary cohomology of {$B_G U(1)$}},
	year = {2017},
	note = {\url{https://arxiv.org/abs/1708.06009}}
}

@book{costenoble_equivariant_2016,
    AUTHOR = {Costenoble, Steven R. and Waner, Stefan},
     TITLE = {Equivariant ordinary homology and cohomology},
    SERIES = {Lecture Notes in Mathematics},
    VOLUME = {2178},
 PUBLISHER = {Springer, Cham},
      YEAR = {2016},
     PAGES = {xiv+294},
      ISBN = {978-3-319-50447-6; 978-3-319-50448-3},
   MRCLASS = {55N91 (55N25 55P42 55P91 55R91 57R91)},
  MRNUMBER = {3585352},
MRREVIEWER = {J.\ P. C. Greenlees},
       DOI = {10.1007/978-3-319-50448-3},
       URL = {https://doi.org/10.1007/978-3-319-50448-3},
       note = {\url{https://arxiv.org/abs/math/0310237}}
}

@article{costenoble_equivariant_1992,
    AUTHOR = {Costenoble, S. R. and Waner, S.},
     TITLE = {The equivariant {T}hom isomorphism theorem},
   JOURNAL = {Pacific J. Math.},
  FJOURNAL = {Pacific Journal of Mathematics},
    VOLUME = {152},
      YEAR = {1992},
    NUMBER = {1},
     PAGES = {21--39},
      ISSN = {0030-8730,1945-5844},
   MRCLASS = {55N91 (55N25 55R91 57R91)},
  MRNUMBER = {1139971},
MRREVIEWER = {L.\ Gaunce\ Lewis, Jr.},
       URL = {http://projecteuclid.org/euclid.pjm/1102636484},
}

@article{costenoble_equivariant_2025a,
	author = {Steven R. Costenoble and Thomas Hudson},
	title = {The {$C_2$}-equivariant ordinary cohomology of complex quadrics {I}: The antisymmetric case},
	year = {2025},
	note = {\url{https://arxiv.org/abs/2510.26361}}
}

@article{costenoble_equivariant_2025b,
	author = {Steven R. Costenoble and Thomas Hudson},
	title = {The {$C_2$}-equivariant ordinary cohomology of complex quadrics {II}: The symmetric case},
	year = {2025},
	note = {\url{https://arxiv.org/abs/2511.14251}}
}

@article{costenoble_equivariant_2026,
    AUTHOR = {Costenoble, Steven R. and Hudson, Thomas},
     TITLE = {The {$C_2$}-equivariant ordinary cohomology of {$BT^2$}},
   JOURNAL = {J. Topol.},
  FJOURNAL = {Journal of Topology},
    VOLUME = {19},
      YEAR = {2026},
    NUMBER = {2},
     PAGES = {Paper No. e70084, 53},
      ISSN = {1753-8416,1753-8424},
   MRCLASS = {55N91 (14M15 55R12 55R40 55R91 57R20)},
  MRNUMBER = {5086795},
       DOI = {10.1112/topo.70084},
       URL = {https://doi.org/10.1112/topo.70084},
       note = {\url{https://arxiv.org/abs/2411.06470}}
}

@article{costenoble_equivariant_2024a,
	author = {Steven R. Costenoble and Thomas Hudson},
	title = {The {$C_2$}-equivariant ordinary cohomology of {$BU(2)$}},
	year = {2024},
	note = {\url{https://arxiv.org/abs/2411.06411}}
}

@article{HKT19,
   title={Anomaly matching in the symmetry broken phase: Domain walls, {CPT}, and the {Smith} isomorphism},
   volume={8},
   ISSN={2542-4653},
   url={http://dx.doi.org/10.21468/SciPostPhys.8.4.062},
   DOI={10.21468/scipostphys.8.4.062},
   number={4},
   journal={SciPost Physics},
   publisher={Stichting SciPost},
   author={Hason, Itamar and Komargodski, Zohar and Thorngren, Ryan},
   year={2020},
   month={Apr},
	note = "\url{https://arxiv.org/abs/1910.14039}"
}

@article{WWZ19,
    AUTHOR = {Wan, Zheyan and Wang, Juven and Zheng, Yunqin},
     TITLE = {Higher anomalies, higher symmetries, and cobordisms {II}:
              {L}orentz symmetry extension and enriched bosonic/fermionic
              quantum gauge theory},
   JOURNAL = {Ann. Math. Sci. Appl.},
  FJOURNAL = {Annals of Mathematical Sciences and Applications},
    VOLUME = {5},
      YEAR = {2020},
    NUMBER = {2},
     PAGES = {171--257},
      ISSN = {2380-288X},
        note = {\url{https://arxiv.org/abs/1912.13504}}
}

@article{milnor_steenrod_1958,
    AUTHOR = {Milnor, John},
     TITLE = {The {S}teenrod algebra and its dual},
   JOURNAL = {Ann. of Math. (2)},
  FJOURNAL = {Annals of Mathematics. Second Series},
    VOLUME = {67},
      YEAR = {1958},
     PAGES = {150--171},
      ISSN = {0003-486X},
   MRCLASS = {55.00 (18.00)},
  MRNUMBER = {99653},
MRREVIEWER = {G.\ Hirsch},
       DOI = {10.2307/1969932},
       URL = {https://doi.org/10.2307/1969932},
}

@phdthesis{yang_thesis,
	title = {Categorical dynamics on Stable Module Categories},
	school = {Harvard University},
	author = {Yang, Lucy},
	year = {2022},
    note = {\url{https://arxiv.org/abs/2212.09964}}
}

@article{ricka_local_2016,
	title = {Local study of stable module categories via tensor triangulated geometry},
	url = {http://arxiv.org/abs/1610.03561},
	doi = {10.48550/arXiv.1610.03561},
	urldate = {2026-08-17},
	publisher = {arXiv},
	author = {Ricka, Nicolas},
	month = oct,
	year = {2016},
	note = {\url{https://arxiv.org/abs/1610.03561}},
}

@article{bhattacharya_stable_2017,
	author = {Prasit Bhattacharya and Nicolas Ricka},
	title = {The stable {P}icard group of {$\mathcal{A}(2)$}},
	year = {2017},
	note = {\url{https://arxiv.org/abs/1702.01493}}
}

@article{pan_stable_2026,
    AUTHOR = {Pan, JianZhong and Yan, RuJia},
     TITLE = {The stable {P}icard group of {$\mathcal {A}(n)$}},
   JOURNAL = {J. Homotopy Relat. Struct.},
  FJOURNAL = {Journal of Homotopy and Related Structures},
    VOLUME = {21},
      YEAR = {2026},
    NUMBER = {1},
     PAGES = {23--43},
      ISSN = {2193-8407,1512-2891},
   MRCLASS = {55S10 (16T05)},
  MRNUMBER = {5015478},
       DOI = {10.1007/s40062-025-00387-4},
       URL = {https://doi.org/10.1007/s40062-025-00387-4},
       note = {\url{https://arxiv.org/abs/2212.09985}}
}

@article{gheorghe_picard_2018,
    AUTHOR = {Gheorghe, Bogdan and Isaksen, Daniel C. and Ricka, Nicolas},
     TITLE = {The {P}icard group of motivic {$\mathcal A_{\mathbb C}(1)$}},
   JOURNAL = {J. Homotopy Relat. Struct.},
  FJOURNAL = {Journal of Homotopy and Related Structures},
    VOLUME = {13},
      YEAR = {2018},
    NUMBER = {4},
     PAGES = {847--865},
      ISSN = {2193-8407,1512-2891},
   MRCLASS = {14F42 (14C22 20C05)},
  MRNUMBER = {3870774},
MRREVIEWER = {Benjamin\ Collas},
       DOI = {10.1007/s40062-018-0200-z},
       URL = {https://doi.org/10.1007/s40062-018-0200-z},
       note = {\url{https://arxiv.org/abs/1606.05191}}
}

@article{gao_picard_2024,
    AUTHOR = {Gao, Xu and Li, Ang},
     TITLE = {The stable {P}icard group of finite {A}dams {H}opf algebroids
              with an application to the {$\mathbb{R}$}-motivic {S}teenrod
              subalgebra {$\mathcal{A}^{\mathbb{R}}(1)$}},
   JOURNAL = {J. Pure Appl. Algebra},
  FJOURNAL = {Journal of Pure and Applied Algebra},
    VOLUME = {228},
      YEAR = {2024},
    NUMBER = {11},
     PAGES = {Paper No. 107732, 21},
      ISSN = {0022-4049,1873-1376},
   MRCLASS = {14F42 (14C22 20G05 55P42 55S10)},
  MRNUMBER = {4751191},
       DOI = {10.1016/j.jpaa.2024.107732},
       URL = {https://doi.org/10.1016/j.jpaa.2024.107732},
       note = {\url{https://arxiv.org/abs/2306.12527}}
}

@incollection{hopkins_constructions_1994,
    AUTHOR = {Hopkins, Michael J. and Mahowald, Mark and Sadofsky, Hal},
     TITLE = {Constructions of elements in {P}icard groups},
 BOOKTITLE = {Topology and representation theory ({E}vanston, {IL}, 1992)},
    SERIES = {Contemp. Math.},
    VOLUME = {158},
     PAGES = {89--126},
 PUBLISHER = {Amer. Math. Soc., Providence, RI},
      YEAR = {1994},
      ISBN = {0-8218-5165-9},
   MRCLASS = {55P42 (55P60 55Q10)},
  MRNUMBER = {1263713},
MRREVIEWER = {Dominique\ Arlettaz},
       DOI = {10.1090/conm/158/01454},
       URL = {https://doi.org/10.1090/conm/158/01454},
}

@article{may_picard_2001,
    AUTHOR = {May, J. P.},
     TITLE = {Picard groups, {G}rothendieck rings, and {B}urnside rings of
              categories},
   JOURNAL = {Adv. Math.},
  FJOURNAL = {Advances in Mathematics},
    VOLUME = {163},
      YEAR = {2001},
    NUMBER = {1},
     PAGES = {1--16},
      ISSN = {0001-8708,1090-2082},
   MRCLASS = {18D10 (55P42)},
  MRNUMBER = {1867201},
MRREVIEWER = {Pilar\ C.\ Carrasco},
       DOI = {10.1006/aima.2001.1996},
       note = {\url{https://www.math.uchicago.edu/~may/PAPERS/PicJan01.pdf}}
}

@article{van_der_meer_endotrivial_2021,
	author = {Jeroen van der Meer and Richard Wong},
	title = {Endotrivial modules for cyclic {$p$}-groups and generalized quaternion groups via {G}alois descent},
	year = {2021},
	note = {\url{https://arxiv.org/abs/2107.06308}}
}

@article{conway_perfect_1968,
    AUTHOR = {Conway, J. H.},
     TITLE = {A perfect group of order {$8,315,553,613,086,720,000$} and the
              sporadic simple groups},
   JOURNAL = {Proc. Nat. Acad. Sci. U.S.A.},
  FJOURNAL = {Proceedings of the National Academy of Sciences of the United
              States of America},
    VOLUME = {61},
      YEAR = {1968},
     PAGES = {398--400},
      ISSN = {0027-8424},
   MRCLASS = {20.27},
  MRNUMBER = {237634},
MRREVIEWER = {J.\ E.\ McLaughlin},
       DOI = {10.1073/pnas.61.2.398},
       URL = {https://doi.org/10.1073/pnas.61.2.398},
}

@article{conway_group_1969,
    AUTHOR = {Conway, J. H.},
     TITLE = {A group of order {$8,315,553,613,086,720,000$}},
   JOURNAL = {Bull. London Math. Soc.},
  FJOURNAL = {The Bulletin of the London Mathematical Society},
    VOLUME = {1},
      YEAR = {1969},
     PAGES = {79--88},
      ISSN = {0024-6093,1469-2120},
   MRCLASS = {20.29},
  MRNUMBER = {248216},
MRREVIEWER = {W.\ J.\ Wong},
       DOI = {10.1112/blms/1.1.79},
       URL = {https://doi.org/10.1112/blms/1.1.79},
}

@article{JFT20,
    AUTHOR = {Johnson-Freyd, Theo and Treumann, David},
     TITLE = {{${\rm H}^4({\rm Co}_0; {\bf Z}) = {\bf Z}/24$}},
   JOURNAL = {Int. Math. Res. Not. IMRN},
  FJOURNAL = {International Mathematics Research Notices. IMRN},
      YEAR = {2020},
    NUMBER = {21},
     PAGES = {7873--7907},
      ISSN = {1073-7928,1687-0247},
   MRCLASS = {20D08 (20J06 55R40)},
  MRNUMBER = {4176841},
MRREVIEWER = {David\ Benson},
       DOI = {10.1093/imrn/rny219},
       URL = {https://doi.org/10.1093/imrn/rny219},
       note = {\url{https://arxiv.org/abs/1707.07587}}
}

@misc{Ranicki_8,
        author = {Andrew Ranicki},
        title = {Eight in algebra, topology and mathematical physics},
        note = {\url{https://webhomes.maths.ed.ac.uk/~v1ranick/eight.htm}}
}

@article{albert_art_2026,
	author = {Jan Albert and Masaki Okada and Yuji Tachikawa},
	title = {The Art of Orbifolding by Huge Non-Abelian Groups},
	year = {2026},
	note = {\url{https://arxiv.org/abs/2610.36024}}
}

@article{balmer_geometry_2025,
    AUTHOR = {Balmer, Paul and Gallauer, Martin},
     TITLE = {The geometry of permutation modules},
   JOURNAL = {Invent. Math.},
  FJOURNAL = {Inventiones Mathematicae},
    VOLUME = {241},
      YEAR = {2025},
    NUMBER = {3},
     PAGES = {841--928},
      ISSN = {0020-9910,1432-1297},
   MRCLASS = {20C20 (13A02 18G80 20J06)},
  MRNUMBER = {4946248},
       DOI = {10.1007/s00222-025-01350-z},
       URL = {https://doi.org/10.1007/s00222-025-01350-z},
       note = {\url{https://arxiv.org/abs/2307.04398}}
}

@article{miller_permutation_2025,
	author = {Sam K. Miller},
	title = {Permutation twisted cohomology, remixed},
	year = {2025},
	note = {\url{https://arxiv.org/abs/2509.00954}}
}

@article{gomez_classification_2026,
	author = {Gómez, Juan Omar and Sam K. Miller},
	title = {The classification of integral endotrivial complexes},
	year = {2026},
	note = {\url{https://arxiv.org/abs/2605.31128}}
}

@article{gomez_permutation_2026,
	author = {Gómez, Juan Omar and Sam K. Miller},
	title = {Permutation twisted cohomology, reduced and immersed},
	year = {2026},
	note = {\url{https://arxiv.org/abs/2610.40602}}
}

@book{EKMM97,
    AUTHOR = {Elmendorf, A. D. and Kriz, I. and Mandell, M. A. and May, J.
              P.},
     TITLE = {Rings, modules, and algebras in stable homotopy theory},
    SERIES = {Mathematical Surveys and Monographs},
    VOLUME = {47},
      NOTE = {With an appendix by M. Cole. \url{https://www.math.uchicago.edu/~may/BOOKS/EKMM.pdf}},
 PUBLISHER = {American Mathematical Society, Providence, RI},
      YEAR = {1997},
     PAGES = {xii+249},
      ISBN = {0-8218-0638-6},
   MRCLASS = {55N20 (19D10 19D55 55P42 55T25)},
  MRNUMBER = {1417719},
MRREVIEWER = {Donald M. Davis},
       DOI = {10.1090/surv/047},
       URL = {https://doi.org/10.1090/surv/047},
}

\end{document}